%%%%%%%%%%%%%%%%   Geometry and Topology: 2002-25.tex  %%%%%%%%%%%%%
%%%%        
%%%%         A chain rule in the calculus of homotopy functors   
%%%%             
%%%%                    John R Klein and John Rognes  
%%%%  
%%%%              Published in Volume 6(2002) pages 853-887
%%%%
%%%%                   Publication date 19 December 2002
%%%%
%%%%                        This is a LaTeX file
%%%%
%%%%
%%%%%%%%%%%%%%%%%%                                   %%%%%%%%%%%%%%%%%%%

%%%%%%%%%%%%%%%%%%%%%%%%%%%%%%%%%%%%%%%%%%%%%%%%%%%%%%%%%%%%%%%
%%%%%%%%%%%             gtmacros.tex            %%%%%%%%%%%%%%%
%%%%%%%%%%%             version 1.6             %%%%%%%%%%%%%%% 
%
%                       Colin Rourke   
%
%
%    These macros are recommended for use by authors submitting articles   
%    to Geometry and Topology or to Algebraic and Geometric Topology.  
%    They are intended to be used with plain TeX. Each macro is described 
%    briefly to make it clear how to use it (or to modify it to achieve
%    different results).  If you modify this file then please change its
%    name.  If you modify this file and use the modified file to 
%    format an article for submission to Geometry and Topology or
%    Algebraic and Geometric Topology, then please paste the modified
%    file into your main TeX file.  Do not submit it as a separate file.
%      
%    Instructions on using these macros are also given in  gtmacins.tex  
%    or  gtmacins.ps  or .pdf  available on the gt www pages or by 
%    anonymous ftp from the gt/info/macros directory.
%
%
\magnification=\magstephalf      % Sets default point size to 11pt.
%
%  Basic layout parameters :
%
\vsize=7.5truein                 % Sets text height to 7.5 inches.
\hsize=5.2truein                 % Sets text width to 5.2 inches.
\newskip\stdskip                 % standard vertical space
\stdskip=6pt plus3pt minus3pt    % (slightly more stretchy
\medskipamount=\stdskip          % than the usual \medskip)
\parindent=0pt                   % Paragraphs are non-indented with
\parskip=\stdskip                % a little space between paragraphs. 
\abovedisplayskip=\stdskip       %  Reduces the space
\belowdisplayskip=\stdskip       %  around displays.
\mathsurround=0.75pt             % Gives a little extra space around maths.
\overfullrule=0pt                %  Prevents black boxes
%
%   The following macro is for principal paragraph breaks ie
%   a paragraph break with a slightly larger space :
%
\def\ppar{\par\goodbreak\vskip 8pt plus 4pt minus 4pt}     
%
%  The standard horizontal space for theorems, labels etc :
%
\def\stdspace{\hskip 0.75em plus 0.15em\ignorespaces}
\let\qua\stdspace % useful abbreviation (3/4 of a quad)
%
%%%%%%%%%%%%%%            FONT MACROS            %%%%%%%%%%%%
%
%           The following font macros define the AMS symbol 
%           and Euler-Fraktal fonts for use in text and
%           mathematics with appropriate size changes.
%           They also define two new control sequences  
%           \small  and  \large  (similar to those built
%           into LaTeX) which change the size of all fonts 
%           both in text and maths.  \small  is 10% smaller 
%           than normal and  \large  30% bigger.  The strange
%           size of the \large text fonts (10pt scaled 1315)
%           is because these macros are intended to be used
%           at \magstephalf.  The result is 10pt scaled 1440
%           (\magstep2) which is a standard font size.  If
%           you are borrowing these macros to use them at
%           another basic  \magnification, then you will
%           probably need to change 1315 to 1200 in the eleven
%           places marked ** below.  \large  will then be
%           20% bigger than normal.  Note that at \magstephalf
%           all the fonts come out roughly one point larger
%           than their size as defined in these macros.
%
%           The size-changing macros are based on Knuth's
%           \ninepoint and \eightpoint macros.
%
%
%    The macros are laid out in a way which makes it clear how to
%    add futher fonts (or delete unavailable ones) and how to add
%    further size changes.
%
%    First comes a definition of  \hexnumber  which is needed to
%    refer to font families whose family number is not known :
%
\def\hexnumber#1{\ifcase#1 0\or 1\or 2\or 3\or 4\or 5\or 6\or 7\or 8\or
 9\or A\or B\or C\or D\or E\or F\fi}
%
%     Next we define the AMS symbol-a fonts at 13,10,9,7,6,5 points
%
\font\thirtnmsa=msam10 scaled 1315    %%% **  see note above 
\font\tenmsa=msam10          \font\ninemsa=msam9
\font\sevenmsa=msam7         \font\sixmsa=msam6
\font\fivemsa=msam5
%%%%%%  (add further sizes here if you need them)
%
%    and the standard size family for these fonts
%
\newfam\msafam                  \textfont\msafam=\tenmsa
\scriptfont\msafam=\sevenmsa    \scriptscriptfont\msafam=\fivemsa
\edef\hexa{\hexnumber\msafam}        %  The msa family is  \fam\hexa
\def\msa{\fam\msafam\tenmsa}         %  \msa  switches to this family
%
%    Repeat these steps for the AMS symbol-b fonts
%
\font\thirtnmsb=msbm10 scaled 1315   %%%  ** see note above
\font\tenmsb=msbm10      \font\ninemsb=msbm9
\font\sevenmsb=msbm7     \font\sixmsb=msbm6
\font\fivemsb=msbm5
%%%%%%  (add further sizes here if you need them)
%
\newfam\msbfam                   \textfont\msbfam=\tenmsb       
\scriptfont\msbfam=\sevenmsb     \scriptscriptfont\msbfam=\fivemsb
\edef\hexb{\hexnumber\msbfam}    %  The msb family is \fam\hexb  
\def\msb{\fam\msbfam\tenmsb}     %  \msb switches to this family
%
%        Repeat for the Euler-Fraktal fonts 
%
\font\thirtneufm=eufm10 scaled 1315   %%% **  see note above 
\font\teneufm=eufm10                 \font\nineeufm=eufm9
\font\seveneufm=eufm7                \font\sixeufm=eufm6
\font\fiveeufm=eufm5
%%%%%%  (add further sizes here if you need them)
%
\newfam\eufmfam                    \textfont\eufmfam=\teneufm
\scriptfont\eufmfam=\seveneufm     \scriptscriptfont\eufmfam=\fiveeufm
\edef\hexf{\hexnumber\eufmfam}      % The Euler-Fraktal family is
\def\frak{\fam\eufmfam\teneufm}     % \fam\hexf and \frak switches to this
%
%%%  Add further fonts families here (using the same format) if you need
%    them.  The def of hexnumber is optional (it is only used for
%    \mathchardef 's).
%
%      Now we need to define the standard fonts (which are
%      already defined at 10,7 and 5 point) at 13,9 and 6 point:
%
%      Roman fonts:
\font\thirtnrm=cmr10 scaled 1315    %%%  ** see note above
\font\ninerm=cmr9                   \font\sixrm=cmr6   
%%%%%%  (add further sizes here if you need them)
%
%      Math italic fonts
\font\thirtni=cmmi10 scaled 1315    %%%  ** see note above 
\font\ninei=cmmi9                   \font\sixi=cmmi6  
%%%%%%  (add further sizes here if you need them)
%
%     Symbol fonts
\font\thirtnsy=cmsy10 scaled 1315   %%%  ** see note above
\font\ninesy=cmsy9                  \font\sixsy=cmsy6  
%%%%%%  (add further sizes here if you need them)
%
%     Bold face
\font\thirtnbf=cmbx10 scaled 1315   %%%  ** see note above 
\font\ninebf=cmbx9                  \font\sixbf=cmbx6  
%%%%%%  (add further sizes here if you need them)
%
%     The maths extension font (only defined at text size)
%
\font\thirtnex=cmex10 scaled 1315   %%%  ** see note above
\font\nineex=cmex9                  
%%%%%%  (add further sizes here if you need them)
%
%     Finally three fonts (text italic, slanted and typewriter type)
%     which are also only defined at text size
%
\font\thirtnit=cmti10 scaled 1315  %%%  ** see note above 
\font\nineit=cmti9                  
%%%%%%  (add further sizes here if you need them)
%
\font\thirtnsl=cmsl10 scaled 1315  %%%  ** see note above 
\font\ninesl=cmsl9                  
%%%%%%  (add further sizes here if you need them)
%
\font\thirtntt=cmtt10 scaled 1315  %%%  ** see note above 
\font\ninett=cmtt9                  
%%%%%%  (add further sizes here if you need them)
%
%
%     Now come the two main macros.  What  \small  does is to
%     change all the families of fonts from normal size which is
%     10,7,5  (ie 10pt text, 7pt subscript, 5pt subsubscript)
%     to 9,6,5.  \large  similarly changes to  13,9,7.  To make
%     other size changing macros, choose your three sizes, add
%     font size definitions if necessary and make the obvious changes
%     to one of these macros.  Change  \normalbaselineskip  and
%     \strutbox  dimensions to appropriate sizes as well.  To
%     add further fonts, insert them in each macro, using the
%     AMS fonts as a model.
%      
%
\def\small{%
%
%   redefine the sizes of the roman fonts :
%
\textfont0=\ninerm \scriptfont0=\sixrm \scriptscriptfont0=\fiverm
\def\rm{\fam0\ninerm}%       % ( \rm  sets \ninerm  in text mode
%                            %  and \fam0 in math mode)
%
%   and the math italic fonts :
%
\textfont1=\ninei \scriptfont1=\sixi \scriptscriptfont1=\fivei
%
%   and the symbol fonts :
%
\textfont2=\ninesy \scriptfont2=\sixsy \scriptscriptfont2=\fivesy
%
%   There is only one math extension font :
%
\textfont3=\nineex \scriptfont3=\nineex \scriptscriptfont3=\nineex
%
%   Next the bold font (named rather than numbered) :
%
\textfont\bffam=\ninebf \scriptfont\bffam=\sixbf
\scriptscriptfont\bffam=\fivebf \def\bf{\fam\bffam\ninebf}%
%
%   and the three text-only fonts : 
%
\textfont\itfam=\nineit \def\it{\fam\itfam\nineit}%
\textfont\slfam=\ninesl \def\sl{\fam\slfam\ninesl}%
\textfont\ttfam=\ninett \def\tt{\fam\ttfam\ninett}%
%
%   Now the three new families of AMS fonts :
%
%   AMS symbol-a
%
\textfont\msafam=\ninemsa \scriptfont\msafam=\sixmsa
\scriptscriptfont\msafam=\fivemsa \def\msa{\fam\msafam\ninemsa}%         
%
%   AMS symbol-b
%
\textfont\msbfam=\ninemsb \scriptfont\msbfam=\sixmsb
\scriptscriptfont\msbfam=\fivemsb \def\msb{\fam\msbfam\ninemsb}%         
%
%   Euler-Fraktal font
%
\textfont\eufmfam=\nineeufm  \scriptfont\eufmfam=\sixeufm
\scriptscriptfont\eufmfam=\fiveeufm \def\frak{\fam\eufmfam\nineeufm}%
%
%%%  Add further fonts families here if you need them.
%
%    Reset \normalbaselineskip and \strubox :
%
\normalbaselineskip=11pt%
\setbox\strutbox=\hbox{\vrule height8pt depth3pt width0pt}%
%
%    Set \normalbaselines and \rm (roman) as defaults :
%
\normalbaselines\rm
%
%    Reset some of the basic vertical skips:
%
\stdskip=4pt plus2pt minus2pt    
\medskipamount=\stdskip          
\parskip=\stdskip                
\abovedisplayskip=\stdskip       
\belowdisplayskip=\stdskip       
\def\ppar{\par\goodbreak\vskip 6pt plus 3pt minus 3pt}%     
%
%   And finally reset the size of section heads (see below):
%
\def\section##1{\global\advance\sectionnumber by 1
\vskip-\lastskip\penalty-800\vskip 20pt plus10pt minus5pt 
\egroup{\bf\number\sectionnumber\quad##1}\bgroup\small         
\vskip 6pt plus3pt minus3pt
\nobreak\resultnumber=1}%      % Reset resultnumber at start of section
}    %%%   End of  \small  macro      
%
%   Two useful abbreviations to keep track of \small material:
\def\beginsmall{\bgroup\small}
\let\endsmall\egroup
%
%
%    The \large  macro is similar (comments abbreviated):
%
%
\def\large{%
\textfont0=\thirtnrm \scriptfont0=\ninerm \scriptscriptfont0=\sevenrm
\def\rm{\fam0\thirtnrm}%
\textfont1=\thirtni \scriptfont1=\ninei \scriptscriptfont1=\seveni
\textfont2=\thirtnsy \scriptfont2=\ninesy \scriptscriptfont2=\sevensy
\textfont3=\thirtnex \scriptfont3=\thirtnex \scriptscriptfont3=\thirtnex
\textfont\bffam=\thirtnbf \scriptfont\bffam=\ninebf
\scriptscriptfont\bffam=\sevenbf \def\bf{\fam\bffam\thirtnbf}%
\textfont\itfam=\thirtnit \def\it{\fam\itfam\thirtnit}%
\textfont\slfam=\thirtnsl \def\sl{\fam\slfam\thirtnsl}%
\textfont\ttfam=\thirtntt \def\tt{\fam\ttfam\thirtntt}%
%   AMS symbol-a  :
\textfont\msafam=\thirtnmsa \scriptfont\msafam=\ninemsa
\scriptscriptfont\msafam=\sevenmsa \def\msa{\fam\msafam\thirtnmsa}%         
%   AMS symbol-b  :
\textfont\msbfam=\thirtnmsb \scriptfont\msbfam=\ninemsb
\scriptscriptfont\msbfam=\sevenmsb \def\msb{\fam\msbfam\thirtnmsb}%         
%   Euler-Fraktal font :
\textfont\eufmfam=\thirtneufm  \scriptfont\eufmfam=\nineeufm
\scriptscriptfont\eufmfam=\seveneufm \def\frak{\fam\eufmfam\teneufm}%
%%%% Add further fonts families here if you need them.
%   Reset \normalbaselineskip and \strubox and initialise :
\normalbaselineskip=16pt%
\setbox\strutbox=\hbox{\vrule height11.5pt depth4.5pt width0pt}%
\normalbaselines\rm}%
\let\Large\large   %  for compatibility with latex
%
%   The next two lines define commonly used switches for
%   blackboard bold (\Bbb) and gothic type (\goth).  The   
%   \Bbb  switch is set to work in the same way as in amstex
%   and switches only the next character to blackboard bold.

%
%   To use the new AMS fonts you can either use the control
%   sequences \msa, \msb (alias \Bbb) and \frak (alias \goth) eg :
\def\N{{\frak N}}

   % see the msam font table
%
%   or, more generally, make \mathchardef's (cf Knuth p155) eg :
\mathchardef\plussquare="0\hexa01
\mathchardef\nge="3\hexb0B
\mathchardef\maltesecross="0\hexa7A
\mathchardef\del="0\hexf01
%
%   or you can use the amstex names for all the new symbols by
%   inserting the line  \input amsnames  in your file directly
%   after \input gtmacros. 
%   This presupposes that you have collected a copy of the file
%   amsnames.tex  from the  gt/info/macros  ftp directory.
%
%
%   Finally we need a small capital font (for author(s)) :
%
\font\sc=cmcsc10
%
%%%%%%%%%%%%%%%%%       END OF FONT MACROS     %%%%%%%%%%%%%
%
%
%                 Knuth's \square macro :
%
\def\sqr#1#2{{\vcenter{\vbox{\hrule  height.#2truept
	\hbox{\vrule width.#2truept height#1truept 
	\kern#1truept \vrule width.#2truept}
	\hrule height.#2truept}}}}
\def\sq{\sqr55}    %   A small square for end-of-proofs. 
%                  %   (Define other size squares by varing the
%                  %   the two numbers.)
%
%
%      Style macros for section heads, theorem statements etc :
%   
%
\newcount\sectionnumber            %%%  Allocate registers to take
\newcount\resultnumber             %%%  section and result numbers.
\sectionnumber=0\resultnumber=1    %%%  Set these registers to 0 and 1
%
%   The \section macro produces a \large bold faced section heading
%   numbered to the left.  Pagebreaks are encouraged before the
%   start of the section and discouraged directly after the heading.
%   Typical use  \section{First steps}  with typical result :
%
%    1  First Steps     (set bold and \large)
%
\def\section#1{\global\advance\sectionnumber by 1
\xdef\nextkey{\number\sectionnumber}%      (used by the \key macro)
\vskip-\lastskip\penalty-800\vskip 20pt plus10pt minus5pt 
{\large\bf\number\sectionnumber\quad#1}         
\vskip 8pt plus4pt minus4pt
\nobreak\resultnumber=1}      % Reset resultnumber at start of section
%
%
%
%   Next a macro to set subheadings (like the  \section  macro
%   but without the number, with less space and set in standard size).
%
%   Typical use :  \sh{Example formats}
%
\def\sh#1{\vskip-\lastskip\ppar{\bf #1}\par\nobreak\medskip}         
%
%   The \proc ... \endproc macros ("proclaim") are for setting theorems, 
%   lemmas, conjectures etc with automatic numbering.  Typical use :    
%  
%    \proc{Theorem}Every lemon is yellow.\endproc
%
%   Typical result :
%     
%    Theorem 3.4  Every lemon is yellow.   

%   (with Theorem 3.4 set bold and a \stdspace of space before the 
%   statement set in slanted type).
%
\def\proc#1{\xdef\nextkey{\number\sectionnumber.\number\resultnumber}%
\vskip-\lastskip\ppar\bf%
\noindent#1\ \number\sectionnumber.\number\resultnumber
\stdspace\sl\global\advance\resultnumber by 1\ignorespaces}
 
%
%  The \prf ... \endprf macros are for setting proofs.  The code
%  for \prf includes the code for \endproc, so there is no need to
%  type \endproc if the theorem is followed immediatedly by a proof.
%
                            %  For start of proofs  
\def\qed{\hfill$\sq$\par\goodbreak\rm}   %  For end (or absence) of proofs
                 %  extra vertical space)  
        %  For start of proof with alternative name
              %  \endproof is an alias for \endprf
%
%   Typical uses :    
%  
%    \proc{Theorem}Every lemon is yellow. \qed\endproc
%
%    \proc{Theorem}Every lemon is yellow.
%    \prf Use your eyes. \endprf
%
%    \proc{Theorem}Every lemon is yellow.
%    \proof{Proof of theorem} Use your eyes. \endprf
%
%   The next macro is a variant of the \proc macro.  It has
%   exactly the same result except that it omits the number.
%
%   Typical use :  
%    
%    \proclaim{Conjecture}Some oranges are yellow.\endproc
%
\def\proclaim#1{\vskip-\lastskip\ppar\bf%
\noindent#1\stdspace\sl\ignorespaces} 
\let\endproclaim\endproc
%
%   The next macro is a further variant for remarks, definitions etc.   
%   It omits the number and does not switch on slanted type.  
%  
%   Typical use :
%
%    \rk{Remark}Some lemons are thick-skinned.\endrk
%
\def\rk#1{\vskip-\lastskip\ppar{\bf #1}\stdspace\ignorespaces}                

%
%   The next macro is for numbering equations etc, \label  produces the 
%   correct number  x.y  and advances the resultnumber register
%
%   Typical use :
%
%     $$fx=7\eqno{\bf\label}$$
%
%   result :
%
%                           fx = 7                           3.5
%
\def\label{\xdef\nextkey{\number\sectionnumber.\number\resultnumber}%
\number\sectionnumber.\number\resultnumber
\global\advance\resultnumber by 1}
%
%
%
%   The next macros are to automate external references.  To use them 
%   type \reflist ..... \endreflist near the beginning of your paper, 
%   where  .... is the list of references in alphabetical order 
%   and in  the form  \key{KEY}  reference    where "KEY" is a 
%   string of characters which reminds you of the reference.   
%   Separate  references with a blank line or a \par.   Eg 
%
%     \reflist
%
%     ..... more references ....
%
%     \key{Kn-84} {\bf D Knuth}, {\it The TeXbook}, Addison--Wesley (1984)
%
%     ..... more references ....
%
%     \endreflist
%
%   Then type  \references  where you wish the references to be printed
%   (normally near the end of the paper).  To refer to Knuth type
%   for example    see Knuth [\ref{Kn-84}, page 320]   and the correct
%   numerical reference will be printed.  Edit the \references macro
%   to change the formatting (if desired).
%   There is an alternative \refkey for \key, provided your KEY contains
%   only letters.  The syntax is:
%
%     \reflist
%
%     ..... more references ....
%
%     \refkey\Knuth  {\bf D Knuth}, {\it The TeXbook}, Addison--Wesley (1984)
%
%     ..... more references ....
%
%     \endreflist
%
%   \key{Knuth}  has exactly the same maening as \refkey\Knuth and you
%   can mix the two syntaxes if you want.  But \refkey\Kn-84
%   would not work.  It would set Kn as the KEY and -84 would get printed!
%
\newcount\refnumber              %  Register for reference numbers
\refnumber=1                     %  set initially to 1.
\long\def\reflist#1\endreflist{%
\long\def\thereflist{#1}{\def\refkey##1##2\par{\xdef##1{\number\refnumber}%
\global\advance\refnumber by 1}%
\def\key##1##2\par{\expandafter\xdef%
\csname##1\endcsname{\number\refnumber}%
\global\advance\refnumber by 1}#1\par}}
\long\def\references{%
\penalty-800\vskip-\lastskip\vskip 15pt plus10pt minus5pt 
{\large\bf References}\ppar %`References' is set \large bold with space around.
{\leftskip=25pt\frenchspacing    % The list of references is set 
\small\parskip=3pt plus2pt       % \small  with small spaces between,
\def\refkey##1##2\par{\noindent  % numbers in [,]'s and set just to the
\llap{[##1]\stdspace}\ignorespaces##2\par}         % left of a 25pt margin.
\def\key##1##2\par{\noindent  
\llap{[\ref{##1}]\stdspace}\ignorespaces##2\par}  
\def\,{\thinspace}\thereflist\par}}
%
%   Next a footnote macro (with automatic numbering) which sets the
%   footnote  \small.
%
%   Typical use :
%         ..... are yellow.\fnote{By yellow here we mean Britsh
%    Standard colour BS3320.} 
%
\newcount\footnotenumber         % Register for footnote number
\footnotenumber=1                % set initially to 1
\def\fnote#1{\xdef\nextkey{\number\footnotenumber}%
{\small\ifnum\footnotenumber>9\parindent=14pt%
\else\parindent=10pt\fi\footnote{$^{\number\footnotenumber}$}%
{\hglue-5pt#1}\global\advance\footnotenumber by 1}}
%
%
%   Next macros for handling figures with automatic numbering (using 
%   TeX's \midinsert to float the figure to a suitable place).
%   
%   The \figure ... \endfigure macro centres the figure and adds
%   an automatically numbered label  Figure XX  after it.
%
%   If you have a caption, then type \caption{caption text} 
%   somewhere between \figure and \endfigure.  The macro
%   will then add  Figure XX: caption text  after the figure.
%
%   If you want an unnumbered or uncentred figure, then use TeX's raw 
%       \midinsert Figure instructions \endinsert  
%   and if you want a numbered figure label in the same style then
%   use \caption{caption text} outside of  \figure ... \endfigure.
%
%   If you need just the label Figure XX  outside of  \figure ... \endfigure
%   then type  \figurelabel .
%
\newcount\figurenumber          % register for figure number
\figurenumber=1                 % set initially to 1
\def\caption#1{\xdef\nextkey{\number\figurenumber}%
\cl{\small Figure \number\figurenumber: #1}%
\global\advance\figurenumber by 1}
\def\figurelabel{\xdef\nextkey{\number\figurenumber}%
\cl{\small Figure \number\figurenumber}%
\global\advance\figurenumber by 1}
\long\def\figure#1\endfigure{{\xdef\nextkey{\number\figurenumber}%
\let\captiontext\relax\def\caption##1{\xdef\captiontext{##1}}%
\midinsert\cl{\ignorespaces#1\unskip\unskip\unskip\unskip}\vglue6pt\cl{\small 
Figure \number\figurenumber\ifx\captiontext\relax\else: \captiontext
\fi}\endinsert\global\advance\figurenumber by 1}}
%
%   Macros for self-correcting internal references.
%
%   There are two macros  \key{KEY}  and  \ref{KEY} .
%
%   The \key macro sets up KEY as a key for whatever number is 
%   being referenced and the \ref macro converts the KEY into 
%   that number.  Type \key after a  \section or \proc or 
%   \label or \fnote or \figure or \caption or \figurelabel .
%
%   Example:
%
%       \section{Introduction}\key{intro}
%       \proc{Theorem}\key{MainTh}Lemons are yelloy\endproc
%       Here we follow\fnote{Follow in the sense of Dickens}
%       \key{Dickens-note}the crowd ....  
%
%       In section \ref{intro}
%       we stated theorem \ref{mainTh} and noted (see footnote 
%       \ref{Dickens-note}) ...
%
\def\nextkey{??}   %  initialise \nextkey (which is reset by all the
%                     numbering macros)
%
\def\key#1{\expandafter\xdef\csname #1\endcsname{\nextkey}}
\def\ref#1{\expandafter\ifx\csname #1\endcsname\relax
\immediate\write16{Reference {#1} undefined}??\else
\csname #1\endcsname\fi}
%
%   Note:  If the KEY contains only letters then \KEY has exactly the
%   same meaning as \ref{KEY} so in the example you could have:
%
%       In section \intro\ we ....
%
%   The \key will work at any time after the macro which sets the
%   number, provided no other macro which sets a number has been used. 
%
%   Macros for forward references:
%              =======
%   The \key \ref macros ONLY work for backwards references.  If you  
%   want to use forwards references, then type \useforwardrefs  near
%   the beginning of your file.  The KEY's are then stored in an
%   auxiliary  .ref  file and you then suffer the same disadvantage as
%   when using LaTeX that you must TeX the file twice to get
%   the references correct.
%
%   To use a forward ref type \ref{KEY}.  (You can type the
%   alternative  \KEY  but you'll get an error on first TeX'ing 
%   if the \KEY is not yet defined.) 
%
%   The macro also allows external references to be listed at the end 
%   of the file (if you wish to).  (Indeed they can be typed anywhere
%   before the \references command.)  You can combine the reference list
%   and the \references command by typing the references (using the
%   same syntax as before) between the commands \biblio and \endbiblio 
%   (don't type \references or they'll be printed twice).
%
\newread\gtinfile
\newwrite\gtreffile
\def\useforwardrefs{
\openin\gtinfile\jobname.ref
\ifeof\gtinfile
\closein\gtinfile
\immediate\write16{No file \jobname.ref}
\else
\closein\gtinfile
\input \jobname.ref
\fi
\immediate\openout\gtreffile \jobname.ref
%
%   Adapt \key :
%
\def\key##1{{\def\\{\noexpand}%
\expandafter\xdef\csname ##1\endcsname{\nextkey}%
\immediate\write\gtreffile{\\\expandafter\\\def\\\csname ##1\\\endcsname%
{\nextkey}}}}
%
%  Adapt macros for external references:  
%
\long\def\reflist##1\endreflist{%
\long\def\thereflist{##1}{\def\refkey####1####2\par{\xdef####1{%
\number\refnumber}{\def\\{\noexpand}\immediate\write\gtreffile
{\\\def\\####1{\number\refnumber}}}\global\advance\refnumber by 1}%
\def\key####1####2\par{\expandafter\xdef%
\csname####1\endcsname{\number\refnumber}%
{\def\\{\noexpand}\immediate\write\gtreffile
{\\\expandafter\\\def\\\csname ####1\\\endcsname{\number\refnumber}}}
\global\advance\refnumber by 1}##1\par}}
\long\def\biblio##1\endbiblio{\reflist##1\endreflist\references}%
%
%  Adapt obselete key macros (\numkey, \seckey and \figkey):
%
\def\numkey##1{{\def\\{\noexpand}%
\xdef##1{\number\sectionnumber.\number\resultnumber}
\immediate\write\gtreffile{\\\def\\##1%
{\number\sectionnumber.\number\resultnumber}}}}
\def\seckey##1{{\def\\{\noexpand}\xdef##1{\number\sectionnumber}
\immediate\write\gtreffile{\\\def\\##1{\number\sectionnumber}}}}
\def\figkey##1{\xdef##1{\number\figurenumber}%
{\def\\{\noexpand}\immediate\write\gtreffile%
{\\\def\\##1{\number\figurenumber}}}
\number\figurenumber\global\advance\figurenumber by 1}
}   %  end of \useforwardrefs
%
%
%   The next five macros are obselete and have been superseeded by
%   the general \key macro above.  They are included merely to 
%   maintain backward compatibility for the package:
%
%
\def\figkey#1{\xdef#1{\number\figurenumber}%
\number\figurenumber\global\advance\figurenumber by 1}
\def\fig#1#2\endfig{%
\midinsert\cl{#2}\vglue6pt\cl{\small Figure #1}\endinsert}
\def\newfig{\number\figurenumber\global\advance\figurenumber by 1}
\def\numkey#1{\xdef#1{\number\sectionnumber.\number\resultnumber}}
\def\seckey#1{\xdef#1{\number\sectionnumber}}
%
%   End of obselete macros.
%
%
%   The next macro is a version of the verbatim macro given by Knuth.
%
%   This macro produces a "verbatim" printout of
%   any ASCII string which does not contain the symbol "
%   (TeX files do not usually contain " 's).
%   More precisely, everything between consecutive pairs
%   of " 's is printed verbatim in the typewriter font cmtt.
%   For an explanation of how the macro works, see Knuth pp 420-1.
%
%   There are two switches: \verb (which switches the macro on)
%   and \brev which switches the macro off (the default).  When
%   the macro is switched off the symbol " has its usual 
%   meaning for TeX.  To use the macro, type \verb before use
%   and the use " to switch verbatim on and off.  Be careful
%   not to use " for any other purpose.  There is no need to
%   switch the macro off again unless you need to use " for
%   some other purpose (eg making  \mathchardef 's).  Note 
%   that the macro MUST BE OFF before inputting  amsnames.tex .
%
%   Whether the macro is on or off you can always use the
%   control sequence \dq (double quote) for " e.g.
%   \mathchardef\sum=\dq1350  is perfectly valid.
%   The control sequence \ttq is an abbreviation for
%   {\tt\dq}.  Thus "\ttq" will produce " (in cmtt)
%   inside a verbatim quote.
%
%
   %  define a code for " so it can be used when \verb is on
  %  code for " in cmtt
%
\def\verb{\catcode`\"=\active}       %  The main
\def\brev{\catcode`\"=12}            %  switches.
\brev                                %  Prime switches and
\verb                                %  switch on.
{\obeyspaces\gdef {\ }}              
{\catcode`\`=\active\gdef`{\relax\lq}}
\def"{%
\begingroup\baselineskip=12pt\def\par{\leavevmode\endgraf}%
\tt\obeylines\obeyspaces\parskip=0pt\parindent=0pt%
\catcode`\$=12\catcode`\&=12\catcode`\^=12\catcode`\#=12%
\catcode`\_=12\catcode`\~=12%
\catcode`\{=12\catcode`\}=12\catcode`\%=12\catcode`\\=12%
\catcode`\`=\active\let"\endgroup}
\brev      %   Finally switch the macro off (for safety)
%
%   Macros for itemised lists.   Typical use :
%    
%    \items
%    \item{(i)}Colours must be defined.
%    \item{(ii)}Colour cards may not be cited.
%    \enditems
%
%   Result :
%
%    (i)  Colours must be defined. 
%   (ii)  Colour cards may not be cited.
%
%
\def\items{\par\leftskip = 25pt}           % Start of itemised list         
\def\enditems{\par\leftskip = 0pt}         % end of itemised list   
\def\item#1{\par\leavevmode\llap{#1\stdspace}%
\ignorespaces}                             % labelled item
               % bulleted item.
%
%   The \quote ... \endquote macros are for typesetting quotations :
%

%
%   A few useful abbreviations :
%
    %  Colon with correct spacing for maps.
\def\np{\vfil\eject}         %  Forced page break (new page).
\def\nl{\hfil\break}         %  New line.
\def\cl{\centerline}         %  Centerline
\def\gt{{\mathsurround=0pt\it $\cal G\mskip-2mu$eometry \&\ 
$\cal T\!\!$opology}}        %  The journal title in recommended style
    %  for monographs
\def\agt{{\mathsurround=0pt\it$\cal A\mskip-.7mu$lgebraic \&\ 
$\cal G\mskip-2mu$eometric $\cal T\!\!$opology}}  % AGT
%
%    Finally some macros for automatic title page or header generation.
%    To use them type your header information using the following  
%    example as a guide :
%
%    Note that \\ is used as standard separator (for lines in \title and
%    \address, between authors and between email addresses or URL's)
%    and that \email, \url and \secondaddress are optional.
%

% Example:  \title{A short spoof paper\\with a two-line title}
% =======   \authors{Albert Einstein\\Leonardo da Vinci}
%           \address{IAS\\Princeton}\secondaddress{Renaissance\\Venice}
%           \email{ae@ias.princeton.edu\\ldv@ren.ven.hist}
%           \abstract 
%           A short spoof paper with a very short abstract.
%           \endabstract 
%           \primaryclass{00-01, 00-02}\secondaryclass{68-00, 68-01}
%           \keywords{Short, spoof, paper}
%           \maketitlepage
%
%
%    The title page or header will then be generated automatically.
%
%
%    Define the various ingredients of the title page:
%
\def\title#1{\def\thetitle{#1}}

\def\author#1{\edef\previousauthors{\theauthors}
 \ifx\theauthors\relax\def\theauthors{#1}\else
 \def\theauthors{\previousauthors\par#1}\fi}

        % aliases
%
\def\address#1{\edef\previousaddresses{\theaddress}
 \ifx\theaddress\relax\def\theaddress{#1}\else
 \def\theaddress{\previousaddresses\par\vskip 2pt\par#1}\fi}
                             % alias
\def\secondaddress#1{\edef\previousaddresses{\theaddress}
 \ifx\theaddress\relax\def\theaddress{#1}\else
 \def\theaddress{\previousaddresses\par{\rm and}\par#1}\fi}   

\def\email#1{\edef\previousemails{\theemail}
 \ifx\theemail\relax\def\theemail{#1}\else
 \def\theemail{\previousemails\hskip 0.75em\relax#1}\fi}
  % aliases
\def\secondemail#1{\edef\previousemails{\theemail}
 \ifx\theemail\relax\def\theemail{#1}\else
 \def\theemail{\previousemails\hskip 0.75em{\rm and}\hskip 0.75em
 \relax#1}\fi}
\def\url#1{\edef\previousurls{\theurl}
 \ifx\theurl\relax\def\theurl{#1}\else
 \def\theurl{\previousurls\hskip 0.75em\relax#1}\fi}
      % aliases
\def\secondurl#1{\edef\previousurls{\theurl}
 \ifx\theurl\relax\def\theurl{#1}\else
 \def\theurl{\previousurls\hskip 0.75em{\rm and}\hskip 0.75em
 \relax#1}\fi}
\long\def\abstract#1\endabstract{\long\def\theabstract{#1}}
\def\primaryclass#1{\def\theprimaryclass{#1}}
\let\subjclass\primaryclass                        % alias
\def\secondaryclass#1{\def\thesecondaryclass{#1}}
\def\keywords#1{\def\thekeywords{#1}}
%
%  Set \\ to \par and title page items to \relax to initialise macros :
%
\let\\\par\let\thetitle\relax\let\theshorttitle\relax
\let\theauthors\relax\let\theshortauthors\relax
\let\theaddress\relax\let\theshortaddress\relax
\let\theemail\relax\let\theurl\relax
\let\theabstract\relax\let\theprimaryclass\relax
\let\thesecondaryclass\relax\let\thekeywords\relax
%
%
%
%   Basic title page layout (edit this macro if you
%   wish to adjust the title page layout) :
%
\long\def\maketitlepage{    % start of definition of \maketitlepage

\vglue 0.2truein   % top margin

% title :
%
{\parskip=0pt\leftskip 0pt plus 1fil\def\\{\par\smallskip}{\large
\bf\thetitle}\par\medskip}   

\vglue 0.15truein 

% authors :
%
{\parskip=0pt\leftskip 0pt plus 1fil\def\\{\par}{\sc\theauthors}
\par\medskip}%
 
\vglue 0.1truein 

% address(es) email's and URL's (with switches to detect whether the
% optional items have been used) :
%
{\small\parskip=0pt
{\leftskip 0pt plus 1fil\def\\{\par}{\sl\theaddress}\par}
\ifx\theemail\relax\else  % email address?
\vglue 5pt \def\\{\stdspace{\rm and}\stdspace} 
\cl{Email:\stdspace\tt\theemail}\fi
\ifx\theurl\relax\else    % URL given?
\vglue 5pt \def\\{\stdspace{\rm and}\stdspace} 
\cl{URL:\stdspace\tt\theurl}\fi\par}

\vglue 7pt 

{\bf Abstract}

\vglue 5pt

\theabstract

\vglue 7pt 

{\bf AMS Classification numbers}\quad Primary:\quad \theprimaryclass\par

Secondary:\quad \thesecondaryclass

\vglue 5pt 

{\bf Keywords:}\quad \thekeywords

\np  % page break at the end of the title page

}    % end of definition of \maketitlepage
%
%    % \makeshorttitle (for general preprints) doesn't take a new page
%
\long\def\makeshorttitle{    % start of definition of \makeshorttitle

%\vglue 0.2truein   % top margin

% title :
%
{\parskip=0pt\leftskip 0pt plus 1fil\def\\{\par\smallskip}{\large
\bf\thetitle}\par\medskip}   

\vglue 0.05truein 

% authors :
%
{\parskip=0pt\leftskip 0pt plus 1fil\def\\{\par}{\sc\theauthors}
\par\medskip}%
 
\vglue 0.03truein 

% address(es) email's and URL's (with switches to detect whether the
% optional items have been used) :
%
{\small\parskip=0pt
{\leftskip 0pt plus 1fil\def\\{\par}{\sl\ifx\theshortaddress\relax
\theaddress\else\theshortaddress\fi}\par}
\ifx\theemail\relax\else  % email address?
\vglue 5pt \def\\{\stdspace{\rm and}\stdspace} 
\cl{Email:\stdspace\tt\theemail}\fi
\ifx\theurl\relax\else    % URL given?
\vglue 5pt \def\\{\stdspace{\rm and}\stdspace} 
\cl{URL:\stdspace\tt\theurl}\fi\par}

\vglue 10pt 

% abstract and classification numbers (with switches):

{\small\leftskip 25pt\rightskip 25pt{\bf Abstract}\stdspace\theabstract

{\bf AMS Classification}\stdspace\theprimaryclass
\ifx\thesecondaryclass\relax\else; \thesecondaryclass\fi\par
{\bf Keywords}\stdspace \thekeywords\par}
\vglue 7pt
}    % end of definition of \makeshorttitle
\let\maketitle\makeshorttitle        %% alias
%
%    %%%% \makeagttitle (for AGT) similar to \makeshorttitle but
%         with addresses omitted (they go at the end)
%
%%%% publication info and test defaults:

\def\volumenumber#1{\def\thevolumenumber{#1}}
\def\volumeyear#1{\def\thevolumeyear{#1}}
\def\pagenumbers#1#2{\def\startpage{#1}\def\finishpage{#2}}
\def\published#1{\def\publishdate{#1}}
\def\received#1{\def\receiveddate{#1}}
\def\revised#1{\def\reviseddate{#1}}
\let\reviseddate\relax
%% Defaults for authors to use to check layout
\volumenumber{X}
\volumeyear{20XX}
\pagenumbers{1}{XXX}
\published{XX Xxxember 20XX}

\long\def\makeagttitle{   %%% start of definition of \makeagttitle
\agt\hfill      %   Journal title (top left) 
%   logo placeholder (top right)
\hbox to 60truept{\vbox to 0pt{\vglue -14truept{\bf [Logo here]}\vss}\hss}
\break
{\small Volume \thevolumenumber\ (\thevolumeyear)
\startpage--\finishpage\nl
Published: \publishdate}

\vglue .2truein

% title
{\parskip=0pt\leftskip 0pt plus 1fil\def\\{\par\smallskip}{\large
\bf\thetitle}\par\medskip}   
\vglue 0.05truein 

% authors :
%
{\parskip=0pt\leftskip 0pt plus 1fil\def\\{\par}{\sc\theauthors}
\par\medskip}%
 
\vglue 0.03truein 

%  abstract and classification numbers:

{\small\leftskip 25truept\rightskip 25truept{\bf Abstract}\stdspace\theabstract

{\bf AMS Classification}\stdspace\theprimaryclass
\ifx\thesecondaryclass\relax\else; \thesecondaryclass\fi\par
{\bf Keywords}\stdspace \thekeywords\par}\vglue 7truept

}   %%%% end of definition of \makeagttitle

%%%%% Macro to typeset addresses (typically at the end of the paper)

\def\Addresses{\bigskip
{\small \parskip 0pt \leftskip 0pt \rightskip 0pt plus 1fil \def\\{\par}
\sl\theaddress\par\medskip \rm Email:\stdspace\tt\theemail\par
\ifx\theurl\relax\else\smallskip \rm URL:\stdspace\tt\theurl\par\fi}}

\def\agtart{%   Full mock-up of AGT article style (for authors to test with)
%  get print centerpage:
\hoffset 14truemm
\voffset 31truemm
\font\phead=cmsl9 scaled 950
\font\pnum=cmbx10 scaled 913
\font\pfoot=cmsl9 scaled 950
%  headline and footline
\headline{\vbox to 0pt{\vskip -4.5mm\line{\small\phead\ifnum
\count0=\startpage ISSN numbers are printed here
\hfill {\pnum\folio}\else\ifodd\count0\def\\{ }% 
\ifx\theshorttitle\relax\thetitle\else\theshorttitle\fi\hfill{\pnum\folio}
\else\def\\{ and }{\pnum\folio}\hfill\ifx\theshortauthors\relax\theauthors
\else\theshortauthors\fi\fi\fi}\vss}}
\footline{\vbox to 0pt{\vglue 0mm\line{\small\pfoot\ifnum\count0=\startpage
Copyright declaration is printed here\hfill\else
\agt, Volume \thevolumenumber\ (\thevolumeyear)\hfill\fi}\vss}}
%  force \agttitle
\let\maketitle\makeagttitle\let\makeshorttitle\makeagttitle}

%%%
%%%  This version of  gtoutput.tex  is intended to finish formatting
%%%  papers published in Geometry & Topology and stored in the
%%%  arXiv.   All versions of  gtoutput.tex  are copyright 
%%%  GT Publications and are to be used _only_ for formatting
%%%  the officially published version of G&T papers.
%%%
%%%
%%%                                             Colin Rourke  14.9.2000
%%%
%%%  To create header file  head.xxx  comment out the first \endinput

%  test for latex or plain tex
\def\ifplaintex{\expandafter\ifx\csname documentclass\endcsname\relax}

%  get print centerpage:

\ifplaintex 
\hoffset 14truemm
\voffset 31truemm
\else
\headsep 23pt
\footskip 35pt
\hoffset -4truemm
\voffset 12.5truemm
\fi

%  load pictex if not already loaded :
\expandafter\ifx\csname beginpicture\endcsname\relax
\expandafter\ifx\csname documentclass\endcsname\relax
\input pictex \else\font\fiverm=cmr5
\input prepictex \input pictex \input postpictex \fi\fi

\def\gt{{\mathsurround=0pt\it $\cal G\mskip-2mu$eometry \&\ 
$\cal T\!\!$opology}}        %  journal title in recommended style

\def\gtp{{\mathsurround=0pt\it $\cal G\mskip-2mu$eometry \&\ 
$\cal T\!\!$opology $\cal P\!$ublications}}  % GT publications

%  define the various new ingredients of the title page 

\def\lognumber#1{\def\thelognumber{#1}}
\def\volumenumber#1{\def\thevolumenumber{#1}}
\def\papernumber#1{\def\thepapernumber{#1}}
\def\volumeyear#1{\def\thevolumeyear{#1}}

\def\pagenumbers#1#2{\def\startpage{#1}\def\finishpage{#2}}
\def\published#1{\def\publishdate{#1}}
\def\proposed#1{\def\theproposer{#1}}
\def\seconded#1{\def\theseconders{#1}}
\def\received#1{\def\receiveddate{#1}}
\def\revised#1{\def\reviseddate{#1}}
\def\accepted#1{\def\accepteddate{#1}}

\def\asciiaddress#1{\def\theasciiaddress{#1}}
\def\asciiemail#1{\def\theasciiemail{#1}}
\long\def\asciiabstract#1{\long\def\theasciiabstract{#1}}

%  initialise

\let\\\par\let\thelognumber\relax
\let\thevolumenumber\relax\let\thepapernumber\relax
\let\thevolumeyear\relax\let\thesamplenumber\relax\let\startpage\relax
\let\finishpage\relax\let\publishdate\relax\let\receiveddate\relax
\let\reviseddate\relax\let\accepteddate\relax\let\theasciititle\relax
\let\theasciiauthors\relax\let\theasciiaddress\relax
\let\theasciiabstract\relax
\let\theasciiemail\relax\let\theshortauthors\relax\let\theshorttitle\relax

\long\def\maketitlep{   % start of definition of \maketitlep

\count0=\startpage

\gt\hfill      %   Journal title (top left) 
%    Logo (top right) :
\beginpicture
\setcoordinatesystem units <0.33truein, 0.33truein> point at 2.2 0.9
\setplotsymbol ({$\cal G$})
\plotsymbolspacing=9truept
\circulararc 315 degrees from 0 1 center at 0 0
\setplotsymbol ({$\cal T$})
\circulararc 315 degrees from 1 -1 center at 1 0
\endpicture
%   end of logo
%
\break
{\small\ifx\thesamplenumber\relax % sample?  
Volume \else Sample
\fi\thevolumenumber\ (\thevolumeyear)
\startpage--\finishpage\nl
Published: \publishdate}
\vglue 0.5truein plus 0.4fil minus 0.1truein

% title
{\parskip=0pt\leftskip 0pt plus 1fil\def\\{\par\smallskip}{\ifplaintex\large
\else\Large\fi\bf\thetitle}\par\medskip}   

\vglue 0pt plus 0.1fil 

% authors
{\parskip=0pt\leftskip 0pt plus 1fil\def\\{\par}{\sc\theauthors}
\par\medskip}

\vglue 0pt plus 0.1fil 

%address(es)
{\small\parskip=0pt\let\newline\\
{\leftskip 0pt plus 1fil\def\\{\par}{\sl\theaddress}\par}
\expandafter\ifx\theemail\relax    % email address?
\relax\else\vglue 5pt plus 0.02fil minus 2pt\def\\{\stdspace{\rm 
and}\stdspace} 
\cl{Email:\stdspace\tt\theemail}\fi
\ifx\theurl\relax                  % URL given?
\relax\else\vglue 5pt plus 0.02fil minus 2pt\def\\{\stdspace{\rm 
and}\stdspace}
\cl{URL:\stdspace\tt\theurl}\fi\par}

\vglue 7pt plus 0.3fil minus 3pt

{\bf Abstract}
\vglue 5pt plus 0.1fil minus 2pt

\theabstract

\vglue 7pt plus 0.3fil minus 3pt

{\bf AMS Classification numbers}\quad Primary:\quad \theprimaryclass

Secondary:\quad \thesecondaryclass

\vglue 5pt plus 0.3fil minus 2pt

{\bf Keywords}\quad \thekeywords

\vglue 10pt plus 0.5fil minus 5pt

{\small  Proposed: \theproposer\hfill Received: \receiveddate\nl
Seconded: \theseconders\hfill 
\ifx\reviseddate\relax                         % paper revised?
Accepted: \accepteddate                        % no
\else
Revised: \reviseddate                          % yes
\fi}
\eject
}       %  end of definition of \maketitlep

\let\maketitlepage\maketitlep
\let\maketitle\maketitlepage

%%% Headers and footers

\font\phead=cmsl9 scaled 950
\font\lhead=cmsl9 scaled 1050
\font\pnum=cmbx10 scaled 913
\font\lnum=cmbx10 
\font\pfoot=cmsl9 scaled 950
\font\lfoot=cmsl9 scaled 1050
\ifplaintex
\headline{\vbox to 0pt{\vskip -4.5mm\line{\small\phead\ifnum
\count0=\startpage ISSN 1364-0380 (on line)
1465-3060 (printed) \hfill {\pnum\folio}\else\ifodd\count0\def\\{ }% 
\ifx\theshorttitle\relax\thetitle\else\theshorttitle\fi\hfill{\pnum\folio}
\else\def\\{ and }{\pnum\folio}\hfill\ifx\theshortauthors\relax\theauthors
\else\theshortauthors\fi\fi\fi}\vss}}
\footline{\vbox to 0pt{\vglue 0mm\line{\small\pfoot\ifnum\count0=\startpage
\copyright\ \gtp\hfill\else
\gt, Volume \thevolumenumber\ (\thevolumeyear)\hfill\fi}\vss
}}
\else
\makeatletter
\def\@oddhead{{\small\lhead\ifnum\count0=\startpage ISSN 1364-0380 (on line)
1465-3060 (printed) \hfill {\lnum\number\count0}\else\ifodd\count0
\def\\{ }\ifx\theshorttitle\relax \thetitle \else\theshorttitle\fi\hfill
{\lnum\number\count0}\else\def\\{ and }{\lnum\number\count0}
\hfill\ifx\theshortauthors\relax 
\theauthors\else\theshortauthors\fi\fi\fi}}\def\@evenhead{\@oddhead}
\def\@oddfoot{\small\lfoot\ifnum\count0=\startpage\copyright\ \gtp\hfill\else
\gt, Volume \thevolumenumber\ (\thevolumeyear)\hfill\fi}
\def\@evenfoot{\@oddfoot}
\makeatother
\fi

   %%%comment out to create xxx header file

\newwrite\gtoutfile
\long\gdef\makeheadfile{  %%% start of definition of \makeheadfile
{\def\\{, }\def\s{ }
\immediate\openout\gtoutfile head.xxx
\immediate\write\gtoutfile{To: math@arxiv.org}
\immediate\write\gtoutfile{Subject: put or rep NNNNN:pppp}
\immediate\write\gtoutfile{--text follows this line--}
\immediate\write\gtoutfile{Proxy-for: \ifx\theasciiauthors\relax
\theauthors\else\theasciiauthors\fi\s<\ifx\theasciiemail\relax\theemail\else\theasciiemail\fi>}
\immediate\write\gtoutfile{\noexpand\\}
\immediate\write\gtoutfile{Authors: \ifx\theasciiauthors\relax
\theauthors\else\theasciiauthors\fi}
{\def\\{ }\immediate\write\gtoutfile{Title: \ifx\theasciititle\relax
\thetitle\else\theasciititle\fi}}
\immediate\write\gtoutfile{Subj-class: GT or SG or MG etc}
\immediate\write\gtoutfile{MSC-class: \theprimaryclass\ifx\thesecondaryclass\relax\else, \thesecondaryclass\fi}
\immediate\write\gtoutfile{Journal-ref: Geom. Topol. \thevolumenumber
(\thevolumeyear) \startpage-\finishpage}
\immediate\write\gtoutfile{Comments: Published by Geometry and Topology at}
\immediate\write\gtoutfile{\s\s http://www.maths.warwick.ac.uk/gt/GTVol\thevolumenumber/paper\thepapernumber.abs.html}
\immediate\write\gtoutfile{\noexpand\\}
\immediate\write\gtoutfile{}
\ifx\theasciiabstract\relax
\immediate\write\gtoutfile{\theabstract}\else
\immediate\write\gtoutfile{\theasciiabstract}\fi
\immediate\write\gtoutfile{}
\immediate\write\gtoutfile{\noexpand\\}
\immediate\write\gtoutfile{}
\immediate\closeout\gtoutfile}}  %%% end of definition of \makeheadfile

\def\maketitlepage{\maketitlep\makeheadfile}
\let\maketitle\maketitlepage

\lognumber{21}

\volumenumber{6}
\papernumber{25} 
\volumeyear{2002}
\pagenumbers{853}{887} 
\received{19 June 1997}
\revised{21 July 2002}
\accepted{19 December 2002}
\published{19 December 2002}
\proposed{Ralph Cohen}
\seconded{Gunnar Carlsson, Thomas Goodwillie}

\input amsnames
\input amstex
\let\cal\Cal     
\catcode`\@=12   
\let\\\par
\def\topmatter{\relax}
\def\endtopmatter{\maketitlepage}
\let\gttitle\title
\def\title#1\endtitle{\gttitle{#1}}
\let\gtauthor\author
\def\author#1\endauthor{\gtauthor{#1}}
\let\gtaddress\address
\def\address#1\endaddress{\gtaddress{#1}}
\let\gtemail\email
\def\email#1\endemail{\gtemail{#1}}
\def\subjclass#1\endsubjclass{\primaryclass{#1}}
\let\gtkeywords\keywords
\def\keywords#1\endkeywords{\gtkeywords{#1}}
\def\heading#1\endheading{{\def\S##1{\relax}\def\\{\relax\ignorespaces}
    \section{#1}}}
\def\head#1\endhead{\heading#1\endheading}

\def\subhead#1\endsubhead{\sh{#1}}
\def\subsubhead#1\endsubsubhead{\sh{#1}}
\def\specialhead#1\endspecialhead{\sh{#1}}
\def\demo#1{\rk{#1}\ignorespaces}
\def\enddemo{\ppar}
\let\remark\demo
\def\endremark{}
\let\definition\demo
\def\enddefinition{\ppar}
\let\example\demo
\def\endexample{\ppar}
\def\qed{\ifmmode\quad\sq\else\hbox{}\hfill$\sq$\par\goodbreak\rm\fi}  
\def\proclaim#1{\rk{#1}\sl\ignorespaces}
\def\endproclaim{\rm\ppar}
\def\cite#1{[#1]}
\newcount\itemnumber
\def\roster{\items\itemnumber=1}
\def\endroster{\enditems}
\let\itemold\item
\def\item{\itemold{{\rm(\number\itemnumber)}}%
\global\advance\itemnumber by 1\ignorespaces}
\def\S{section~\ignorespaces}  %%  expand \S to "section"
\def\date#1\enddate{\relax}
\def\thanks#1\endthanks{\relax}   %%%  Move acknowledgements "manually"
\def\dedicatory#1\enddedicatory{\relax}  %%% to the end of the intro.
  % in some versions of amstex but not all
\let\footnote\plainfootnote

%%   Adapt the amstex reference list macros 
\def\Refs{\ppar{\large\bf References}\ppar\bgroup\leftskip=25pt
\frenchspacing\parskip=3pt plus2pt\small}       
\def\endRefs{\egroup}
\def\widestnumber#1#2{\relax}
\def\endrefitem{}
\def\refdef#1#2#3{\def#1{\leavevmode\unskip\endrefitem#2\def\endrefitem{#3}}}
\def\ref{\par}
\def\endref{\endrefitem\par\def\endrefitem{}}
\refdef\key{\noindent\llap\bgroup[}{]\ \ \egroup}
\refdef\no{\noindent\llap\bgroup[}{]\ \ \egroup}
\refdef\by{\bf}{\rm, }
\refdef\manyby{\bf}{\rm, }
\refdef\paper{\it}{\rm, }
\refdef\book{\it}{\rm, }
\refdef\jour{}{ }
\refdef\vol{}{ }
\refdef\yr{$(}{)$ }
\refdef\ed{(}{ Editor) }
\refdef\publ{}{ }
\refdef\inbook{from: ``}{'', }
\refdef\pages{}{ }
\refdef\page{}{ }
\refdef\paperinfo{}{ }
\refdef\bookinfo{}{ }
\refdef\publaddr{}{ }
\refdef\eds{(}{ Editors)}
\refdef\bysame{\hbox to 3 em{\hrulefill}\thinspace,}{ }
\refdef\toappear{(to appear)}{ }
\refdef\issue{no.\ }{ }

%%      Macros to change refs to numbers
\newcount\refnumber\refnumber=1
\def\refkey#1{\expandafter\xdef\csname cite#1\endcsname{\number\refnumber}%
\global\advance\refnumber by 1}
\def\cite#1{[\csname cite#1\endcsname]}
\def\Cite#1{\csname cite#1\endcsname}  %% unbracketed \cite 
\def\key#1{\noindent\llap{[\csname cite#1\endcsname]\ \ }}
\refkey{Ar99}
\refkey{BHM93}
\refkey{BF78}
\refkey{BK72}
\refkey{Fa95}
\refkey{Go90}
\refkey{Go92}
\refkey{Ha02}
\refkey{He92}
\refkey{Kl95}
\refkey{Ly98}
\refkey{Ma71}
\refkey{Sch01}
\refkey{Wa85}
\refkey{Wa96}

\define\In{\operatorname{in}}
\define\Map{\operatorname{Map}}
\define\OO{\Cal O}
\define\SP{\Cal Sp}
\define\SS{\Cal S}
\define\THH{\operatorname{THH}}
\define\Sing{\operatorname{Sing}}
\define\TT{\Cal T}
\define\UU{\Cal U}
\define\colim{\operatornamewithlimits{colim}}
\define\diag{\operatorname{diag}}
\define\hocolim{\operatornamewithlimits{hocolim}}
\define\hofib{\operatorname{hofib}}
\define\holim{\operatornamewithlimits{holim}}
\define\id{\operatorname{id}}
\define\pr{\operatorname{pr}}

\def\:{\kern.1ex\colon\thinspace}

\topmatter
\title A chain rule in the calculus of homotopy functors \endtitle
\author John R Klein\\John Rognes \endauthor
\address Department of Mathematics, Wayne State University\\Detroit,
Michigan 48202, USA\endaddress
\secondaddress{Department of Mathematics, University of 
Oslo\\N--0316 Oslo, Norway}
\asciiaddress{Department of Mathematics, Wayne State University\\Detroit,
Michigan 48202, USA\\and\\Department of Mathematics, University of 
Oslo\\N--0316 Oslo, Norway}

\email klein@math.wayne.edu\endemail
\secondemail{rognes@math.uio.no}
\asciiemail{klein@math.wayne.edu, rognes@math.uio.no}

\abstract
We formulate and prove a chain rule for the {\sl derivative}, in the
sense of Goodwillie, of compositions of weak homotopy functors from
simplicial sets to simplicial sets.  The derivative spectrum $\partial
F(X)$ of such a functor $F$ at a simplicial set $X$ can be equipped with a
right action by the loop group of its domain $X$, and a free left action
by the loop group of its codomain $Y = F(X)$.  The derivative spectrum
$\partial (E \circ F)(X)$ of a composite of such functors is then stably
equivalent to the balanced smash product of the derivatives $\partial
E(Y)$ and $\partial F(X)$, with respect to the two actions of the loop
group of $Y$.  As an application we provide a non-manifold computation
of the derivative of the functor $F(X) = Q(\Map(K, X)_+)$.
\endabstract

\asciiabstract{We formulate and prove a chain rule for the derivative, 
in the sense of Goodwillie, of compositions of weak homotopy functors
from simplicial sets to simplicial sets.  The derivative spectrum
dF(X) of such a functor F at a simplicial set X can be
equipped with a right action by the loop group of its domain X, and
a free left action by the loop group of its codomain Y = F(X).  The
derivative spectrum d(E o F)(X)$ of a composite of such
functors is then stably equivalent to the balanced smash product of
the derivatives dE(Y) and dF(X), with respect to
the two actions of the loop group of Y.  As an application we
provide a non-manifold computation of the derivative of the functor
F(X) = Q(Map(K, X)_+).}

\primaryclass{55P65}
\secondaryclass{55P42, 55P91}
\keywords Homotopy functor, chain rule, Brown representability \endkeywords
\endtopmatter

\catcode`\@=\active

\input xyall
\SelectTips{cm}{}
\newdir{ >}{{}*!/-10pt/@{>}}

\document

\head Introduction \endhead

The calculus of functors was introduced by Goodwillie in \cite{Go90}
as a language to keep track of stable range calculations of certain
geometrically defined homotopy functors, such as stable pseudo-isotopy
theory.  The input for the theory is a homotopy functor
$$
f \: \UU @>>> \TT
$$
from spaces to based spaces.  At an object $X \in \UU$ it is then
possible to associate the ``best excisive approximation'' to $f$ near $X$.
This so-called {\sl linearization\/} of $f$ at $X$ is a functor
$$
P_X f \: \UU/X @>>> \TT
$$
from spaces over $X$ to based spaces, which maps homotopy pushout
squares to homotopy pullback squares.  The associated reduced functor
is called the {\sl differential\/} of $f$ at $X$, and is denoted
by $D_X f$.  Choosing a base point $x \in X$, the composite functor
$$
L \: \TT @>i>> \UU/X @>D_Xf>> \TT
$$
that takes a based space $T$ to $D_X f(X \vee_x T) = \hofib(P_X f(X
\vee_x T) \to P_X f(X))$ is a {\sl linear\/} functor, whose homotopy
groups $L_*(T) = \pi_*(L(T))$ define a generalized homology theory.
Each such homology theory is represented by a spectrum, and the spectrum
associated to this particular homology theory $L_*$ is called the {\sl
derivative} $\partial f(X)$ of $f$ at $(X, x)$.

The goal of this paper is to establish a chain rule for the derivative of
a composite functor.  This is a reasonable goal, since many naturally
occurring functors are composites.  For example, the topological
Hochschild homology $\THH(X)$ of a space $X$ has the homotopy type of
$Q(\Lambda X_+)$, where $\Lambda X = \Map(S^1, X)$ is the free loop
space (see [\Cite{BHM93}, 3.7]).  We can view this as the composite of the
two functors $f(X) = \Lambda X$ and $e(Y) = Q(Y_+) = \colim_n \Omega^n
\Sigma^n(Y_+)$.

In order to even state a chain rule, some modification has to be made
to the above set-up.  In particular, we will relax the condition that
the functor $f$ takes values in based spaces, considering instead weak
homotopy functors
$$
f \: \UU @>>> \UU
$$
from spaces to spaces.  (All our spaces will be compactly generated.)
Then for any space $X$ we let $Y = f(X)$, and choose base points $x
\in X$ and $y \in Y$.  We then study the derivative $\partial^x_y f(X)$
of $f$ at $X$, with respect to the base points $x$ and $y$.  Of course,
if $f(X)$ naturally comes equipped with a base point, then we may take
that point as $y$.

Thus consider functors $e, f \: \UU \to \UU$, with composite $e \circ
f \: \UU \to \UU$.  Let $X$ be a space, and set $Y = f(X)$, $Z = e(Y)$.
Choose base points $x \in X$, $y \in Y$ and $z \in Z$.  Suppose that $f$
and $e$ are bounded below, stably excisive functors (section~3), that $e$
satisfies the colimit axiom (section~3), and that $Y$ is path connected.
Let $\Omega_y(Y)$ denote the geometric realization of the Kan loop group
(section~8) of the total singular simplicial set of $Y$.  This is a
topological group, weakly homotopy equivalent to the usual loop space
of $(Y, y)$.  (More precisely, $\Omega_y(Y)$ is a group object in the
category of compactly generated topological spaces.)  Then it turns out
that, by choosing the models right (section~9, see also remark~12.4),
the derivative $\partial^x_y f(X)$ admits a left $\Omega_y(Y)$--action and
the derivative $\partial^y_z e(Y)$ admits a right $\Omega_y(Y)$--action.
It thus makes sense to form the homotopy orbit spectrum for the diagonal
$\Omega_y(Y)$--action on the smash product of spectra $\partial^y_z e(Y)
\wedge \partial^x_y f(X)$.

\proclaim{Theorem 1.1}{\rm(Chain Rule)}\qua
Let $e, f \: \UU \to \UU$ be bounded below, stably excisive functors,
with $Y = f(X)$ and $Z = e(Y)$, and suppose that $e$ satisfies the colimit
axiom.  Suppose that $Y$ is path connected, and choose base points $x
\in X$, $y \in Y$ and $z \in Z$.  Then the composite $e \circ f$ is
bounded below and stably excisive, and its derivative spectrum at $X$
with respect to $x$ and $z$ is described by a stable equivalence
$$
\partial^x_z (e \circ f)(X) \simeq
\partial^y_z e(Y) \wedge_{h\Omega_y(Y)} \partial^x_y f(X) \,.
$$
The subscript $h\Omega_y(Y)$ denotes homotopy orbits with respect to
the diagonal action of the topological group $\Omega_y(Y)$.
\endproclaim

This is theorem~12.3 specialized to the case when $Y$ is path connected.

If $X_\alpha$ is the path component of $x$ in $X$, and $Z_\gamma$
is the path component of $z$ in $Z$, then the topological group
$\Omega_x(X_\alpha)$ acts on $\partial^x_z (e \circ f)(X)$
and $\partial^x_y f(X)$ from the right, the topological group
$\Omega_z(Z_\gamma)$ acts on $\partial^x_z (e \circ f)(X)$ and
$\partial^y_z e(Y)$ from the left, and the chain rule gives a stable
equivalence of spectra with left $\Omega_z(Z_\gamma)$--action and right
$\Omega_x(X_\alpha)$--action.

It is technically easier to discuss these group actions on spectra that
are formed in the category $\SS_*$ of based simplicial sets than for
spectra formed in $\TT$.  The reason is that the definition of the right
action by $\Omega_y(Y)$ on $\partial^y_z e(Y)$ basically requires $e$
to be a continuous functor.  It is awkward to achieve continuity from a
weak homotopy functor in the topological context.  However, for functors
between simplicial sets it is easy to promote a weak homotopy functor to
a simplicial functor, which suffices to define the right action in the
simplicial context.  See definition~9.6.  We therefore choose to develop
the whole theory for weak homotopy functors $F \: \SS \to \SS$ from the
category $\SS$ of simplicial sets to itself, rather than for functors
$f \: \UU \to \UU$.  In this case, the chain rule appears as theorem~11.4.

It is also possible to start with functors $\Phi \: \SS/X \to \SS$ and
$\Psi \: \SS/Y \to \SS$, with $X$ a simplicial set and $Y = \Phi(X)$,
that may or may not factor through the forgetful functors $u \: \SS/X
\to \SS$ and $u \: \SS/Y \to \SS$, respectively.  The latter is the most
convenient general framework, and the body of the paper is written in
this context.  Thus theorem~11.3 is really our main theorem, from which
the other forms of the chain rule are easily deduced.

The contents of the paper are as follows.  In section~2 we define the
categories of simplicial sets and spectra that we shall work with,
and fix terminology like ``bounded below'' and ``stably excisive''
in section~3.  Then in section~4 we start with a stably excisive weak
homotopy functor $\Phi \: \SS/X \to \SS$ and construct its ``best excisive
approximation'' $P'\Phi$, adapting [\Cite{Go90}, \S1].  Some modification
is needed, since we want $P'\Phi$ to take values in $\SS/Y$ in order to
be able to compose with $\Psi$.  In section~5 we recall the Goodwillie
derivative $\partial\Phi(X)$.  In section~6 we show that in order to
prove a chain rule expressing $\partial(\Psi \circ \Phi)(X)$ in terms of
$\partial\Psi(Y)$ and $\partial\Phi(X)$ we may replace $\Phi$ and $\Psi$
by their respective best excisive approximations (see proposition~6.2).
This leads us to study linear functors $\Phi \: R(X) \to R(Y)$ and $\Psi
\: R(Y) \to R(Z)$, where $R(X)$ is the category of retractive simplicial
sets over $X$, and $Z = \Psi(Y)$.  In section~7 we reduce further to
the case where $X$, $Y$ and $Z$ are all connected.

When $(Y, y)$ is based and connected, there is a natural equivalence $R(Y)
\simeq R(*, \Omega_y(Y))$ (see proposition~8.1), where $R(*, \Omega_y(Y))$
is the category of based, free $\Omega_y(Y)$--simplicial sets, which we
study in section~8.  We are thus led to study linear functors $\phi \:
R(*, \Omega_x(X)) \to R(*, \Omega_y(Y))$ and $\psi \: R(*, \Omega_y(Y))
\to R(*, \Omega_z(Z))$, and their composite $\psi \circ \phi$.
The Goodwillie derivative $\partial\phi$ of such a functor $\phi$ has a
natural free left $\Omega_y(Y)$--action.  Using a canonical enrichment
of $\phi$ to a simplicial functor $\check\phi$, we show in section~9
that $\partial\phi$ also has a natural right $\Omega_x(X)$--action
(see proposition~9.6).  Then, in section~10 we establish a version of
the {\sl Brown--Whitehead representability theorem\/} (see [\Cite{Go90},
1.3]) that represents a linear functor like $\phi$ above in terms of its
Goodwillie derivative $\partial\phi$, equipped with these left and right
actions, under the assumption that $\phi$ satisfies a ``colimit axiom''.
See propositions~10.1 and~10.4.  In section~11 we bring these structures
and representations together, to prove the chain rule for bounded below,
excisive functors $\Phi$ and $\Psi$ in proposition~11.1, and for bounded
below, stably excisive functors $\Phi$ and $\Psi$ in theorem~11.3.
The translation to functors to and from topological spaces goes via the
usual equivalence $\SS \simeq \UU$, and is found in section~12.

We give a list of examples in section~13, including a purely
homotopy-theoretic derivation in example~13.4 of the ``stable homotopy of
mapping spaces'' functor $X \mapsto Q(\Map(K, X)_+)$, which was previously
investigated in [\Cite{Go90}, \S2], \cite{He92}, and \cite{Ar99} using
manifold or configuration space techniques.  Our answer apparently takes
a different form from that given in the cited papers, but in \cite{Kl95}
the first author shows that the two answers are indeed equivalent.

The paper is written using fairly strong explicit hypotheses on the
functors, such as being bounded below and stably excisive, in line with
the style of \cite{Go90} and \cite{Go92}.  Yet, many of the functors one
typically considers satisfy these hypotheses.  Our main technical reason
for doing so occurs at the end of the proof of proposition~11.1, where we
wish to ensure that one functor respects certain stable equivalences of
spectra arising from another functor.  A side effect is that all proofs
become explicit, appealing directly to homotopy excision rather than
to closed model category theory.  Conceivably some of these conditions
could be relaxed by reference to the framework of simplicial functors,
as in \cite{Ly98}, but the work leading to the present paper precedes
that preprint.  Likewise, the present work can be incorporated into the
more general language of pointed simplicial algebraic theories, as in
\cite{Sch01}.  The second author's Master student H{.} Fausk \cite{Fa95}
proved a version of the chain rule in the special case when $Y = f(X)$
is contractible.

This paper was first written in 1995, following a visit by the first
author to Oslo as a guest of the Norwegian Academy of Sciences.
The authors thank the referees for their detailed and constructive
feedback, and apologize for the long delay in finalizing the manuscript.

\head Categories of simplicial sets \endhead

Let $\SS$ be the category of simplicial sets, and let $X$ be a (fixed)
simplicial set.  The category $\SS/X$ of {\sl simplicial sets over $X$}
has objects the simplicial sets $X'$ equipped with a map $X' \to X$, and
morphisms the maps $X' \to X''$ commuting with the structure map to $X$.
The category $\SS/X$ has the identity map $X \to X$ as a terminal object.

The category $R(X)$ of {\sl retractive simplicial sets over $X$} has
objects the simplicial sets $X'$ with maps $r \: X' \to X$ and $s \:
X \to X'$ such that $rs \: X \to X$ is the identity map, and morphisms
the maps $X' \to X''$ that commute with both structure maps $r$ and $s$.
The category $R(X)$ has the identity map $X \to X$ as an initial and
terminal object.  We briefly denote this base point object by $X$.
In the case $X = *$ (a one-point simplicial set), $R(*)$ is isomorphic
to the category $\SS_*$ of based simplicial sets.

Let $G$ be a simplicial group.  A $G$--simplicial set $W$ is a simplicial
set with an action $G \times W \to W$.  For any $G$--simplicial set
$W$ let $R(W, G)$ be the category of relatively free, retractive
$G$--simplicial sets over $W$.  It has objects $(W', r, s)$, where $W'$
is a $G$--simplicial set, $r \: W' \to W$ and $s \: W \to W'$ are maps
of $G$--simplicial sets, $rs \: W \to W$ is the identity map, and $W'$
may be obtained from $W$ by {\sl attaching free $G$--cells}, ie,
by repeated pushouts along the inclusions $G \times \partial\Delta^n
\subset G \times \Delta^n$.  When $G = 1$ is the trivial group, $R(W, 1)
= R(W)$ as before.  When $W = *$, the objects of $R(*, G)$ are precisely
the based, free $G$--simplicial sets.  (Cf [\Cite{Wa85}, page 378].)

Let $u \: \SS/X \to \SS$, $v \: R(X) \to \SS/X$ and $w \: R(W, G) \to
R(W)$ be the obvious forgetful functors.

Consider any functor $\Phi \: \SS/X \to \SS$.  Let $Y = \Phi(X)$ be
its value at the terminal object $X$ (equipped with the identity map
$X \to X$).  Then there is a canonical lift of $\Phi \: \SS/X \to \SS$
over $u \: \SS/Y \to \SS$ to a functor $\SS/X \to \SS/Y$, which we also
denote by~$\Phi$.  Furthermore, there is a canonical lift of $\Phi \circ
v \: R(X) \to \SS/Y$ over $v \: R(Y) \to \SS/Y$ to a functor $R(X) \to
R(Y)$, which we again denote by~$\Phi$.  The latter functor takes the
chosen initial and terminal object $X$ of $R(X)$ to the chosen initial and
terminal object $Y$ of $R(Y)$.  Such functors are called {\sl pointed}.

Functors $\Phi \: \SS/X \to \SS$ sometimes arise from functors $F \:
\SS \to \SS$ as composites $\Phi = F \circ u$, but will in general depend
on the structure map to $X$.  We have a commutative diagram:
$$
\xymatrix@!C{
R(X) \ar[r]^v \ar[d]^{\Phi} & {\SS/X} \ar[r]^u \ar[d]^{\Phi} \ar[dr]^{\Phi}
& {\SS} \ar@{-->}[d]^F \\
R(Y) \ar[r]^v & {\SS/Y} \ar[r]^u & {\SS}
}
\tag 2.1
$$

In this paper, a {\sl spectrum} $\bold L$ is a sequence $\{n \mapsto
L_n\}$ of based simplicial sets $L_n$, and based structure maps $\Sigma
L_n = L_n \wedge S^1 \to L_{n+1}$ for $n\ge0$, as in [\Cite{BF78}, 2.1].
Here it will be convenient to interpret $S^1$ as $\Delta^1 \cup_{\partial
\Delta^1} \Delta^1$, rather than as $\Delta^1/\partial \Delta^1$.  To be
definite, we take the $0$-th vertex of $\partial\Delta^1$ as the base
point of $S^1$.  Let $S^n = S^1 \wedge \dots \wedge S^1$ (with $n$ copies
of $S^1$), and let $CS^n = S^n \wedge \Delta^1$ be the cone on $S^n$.
We write $\SP$ for the category of spectra.

A map of spectra $f \: \bold L \to \bold M$ is a {\sl strict
equivalence\/} if each map $f_n \: L_n \to M_n$ is a weak equivalence.  It
will be called a {\sl meta-stable equivalence\/} if there exist integers
$c$ and $\kappa$ such that $f_n \: L_n \to M_n$ is $(2n-c)$--connected
for all $n \ge \kappa$ (cf section~3).  And $f$ is a {\sl stable
equivalence\/} if it induces an isomorphism $\pi_*(f) \: \pi_*(\bold L)
\to \pi_*(\bold M)$ on all homotopy groups.  Clearly strict equivalences
are meta-stable, and meta-stable equivalences are stable.

Let $G$ be a simplicial group, as above.  A {\sl spectrum with $G$--action}
$\bold L$ is a sequence $\{n \mapsto L_n\}$ of $G$--simplicial sets with
a $G$--fixed base point, and based $G$--maps $\Sigma L_n = L_n \wedge S^1
\to L_{n+1}$ for $n\ge0$, where $G$ acts trivially on $S^1$.  Let $\SP^G$
be the category of spectra with $G$--action.

A {\sl free $G$--spectrum} $\bold L$ is a sequence $\{n \mapsto L_n\}$
of based, free $G$--simplicial sets, and based $G$--maps $\Sigma L_n = L_n
\wedge S^1 \to L_{n+1}$.  Here $G$ acts trivially on $S^1$.  Let $\SP(G)$
be the category of free $G$--spectra.

There are obvious forgetful functors $\SP(G) \to \SP^G$ and $\SP^G
\to \SP$.  A map of free $G$--spectra, or of spectra with $G$--action,
is said to be a ``strict'', ``meta-stable'' or ``stable equivalence''
if the underlying map of spectra has the corresponding property.
In particular, a stable equivalence of spectra with $G$--action is no
more than a $G$--equivariant map that induces an isomorphism on all
homotopy groups.  This naive notion of stable equivalence permits the
formation of homotopy orbits, but not (strict) fixed-points or orbits.

\head Excision conditions \endhead

A morphism $f \: X_0 \to X_1$ in $\SS$ is {\sl $k$--connected\/} if
for every choice of base point $x \in X_0$ the induced map $\pi_n(f)
\: \pi_n(X_0, x) \to \pi_n(X_1, f(x))$ is injective for $0 \le n < k$
and surjective for $0 \le n \le k$.  (No choice of base point is needed
for $n=0$, taking care of the case when $X_0$ is empty.)  A {\sl weak
equivalence\/} is a map that is $k$--connected for every integer $k$.

Let $c$ and $\kappa$ be integers.  A functor $F \: \SS \to \SS$ is said to
satisfy condition $E_1(c, \kappa)$ if for every $k$--connected map $X_0 \to
X_1$ with $k \ge \kappa$ the map $F(X_0) \to F(X_1)$ is $(k-c)$--connected.
A functor satisfying condition $E_1(c, \kappa)$ for some $c$ and $\kappa$
will be called {\sl bounded below}.  Such a functor necessarily takes
weak equivalences to weak equivalences, ie, is a {\sl weak homotopy
functor}.

We form functorial homotopy limits and homotopy colimits of diagrams of
simplicial sets as in \cite{BK72}.  A commutative square of simplicial
sets
$$
\xymatrix{
X_0 \ar[r] \ar[d] & X_1 \ar[d] \\
X_2 \ar[r] & X_3
}
\tag 3.1
$$
is {\sl $k$--cartesian\/} if the induced map $a \: X_0 \to \holim(X_1
\to X_3 \leftarrow X_2)$ is $k$--connected.  It is {\sl cartesian\/}
if $a$ is a weak equivalence.  The square is {\sl $k$--cocartesian\/}
if the induced map $b \: \hocolim(X_1 \leftarrow X_0 \to X_2) \to X_3$
is $k$--connected.  It is {\sl cocartesian\/} if $b$ is a weak equivalence.
(Cf [\Cite{Go90}, 1.2].)

A functor $F \: \SS \to \SS$ is said to satisfy condition $E_2(c, \kappa)$
if, for every cocartesian square as above for which $X_0 \to X_i$ is
$k_i$--connected and $k_i \ge \kappa$ for $i=1,2$, the resulting square
$$
\xymatrix{
F(X_0) \ar[r] \ar[d] & F(X_1) \ar[d] \\
F(X_2) \ar[r] & F(X_3)
}
$$
is $(k_1+k_2-c)$--cartesian.  The functor $F$ is called {\sl stably
excisive\/} if it satisfies condition $E_2(c, \kappa)$ for some integers
$c$ and $\kappa$.  $F$ is called {\sl excisive\/} if it takes
all cocartesian squares to cartesian squares. (Cf [\Cite{Go90}, 1.8].)

A morphism in one of the categories $\SS/X$, $R(X)$ or $R(W, G)$ is said
to be ``$k$--connected'', or a ``weak equivalence'', if the underlying
morphism in $\SS$ has that property.  Similarly for $k$--cartesian,
cartesian, $k$--cocartesian and cocartesian squares.  The conditions
$E_1(c, \kappa)$, ``bounded below'', ``weak homotopy functor'', $E_2(c,
\kappa)$, ``stably excisive'' and ``excisive'' then also make sense for
functors $\SS/X \to \SS$, $\SS/X \to \SS/Y$, $R(X) \to R(Y)$, $R(*, H)
\to R(*, G)$, etc.

\proclaim{Proposition 3.2}
Let $X$ be a simplicial set, $\Phi \: \SS/X \to \SS$ a functor, $Y =
\Phi(X)$ a simplicial set, and $\Psi \: \SS/Y \to \SS$ a functor.
Suppose that $\Phi$ and $\Psi$ are bounded below and stably excisive.
Then the composite functor $\Psi \circ \Phi \: \SS/X \to \SS$ is also
bounded below and stably excisive.
\endproclaim

\demo{Proof}
Suppose that $\Phi$ and $\Psi$ satisfy $E_1(c, \kappa)$ and $E_2(c,
\kappa)$, where we may assume that $c\ge1$ and $\kappa\ge0$.
We claim that $\Psi \circ \Phi$ satisfies $E_1(2c, \kappa+c)$ and
$E_2(3c+1, \kappa+c)$.  The first claim is clear.  For the second,
consider a cocartesian diagram as in~(3.1), with $X_0 \to X_i$
$k_i$--connected for $i=1,2$, and $k_i \ge \kappa+c$.  Apply $\Phi$
to get a $(k_1+k_2-c)$--cartesian square
$$
\xymatrix{
\Phi(X_0) \ar[r] \ar[d] &
\Phi(X_1) \ar[d] \\
\Phi(X_2) \ar[r] &
\Phi(X_3)
}
\tag 3.3
$$
with $\Phi(X_0) \to \Phi(X_i)$ $(k_i-c)$--connected for $i=1,2$.  Let
$$
PO = \hocolim(\Phi(X_1) \leftarrow \Phi(X_0) \rightarrow \Phi(X_2))
$$
be the homotopy pushout in this square.  By homotopy excision
(cf [\Cite{Ha02}, 4.23]) the cocartesian square
$$
\xymatrix{
\Phi(X_0) \ar[r] \ar[d] &
\Phi(X_1) \ar[d] \\
\Phi(X_2) \ar[r] &
PO
}
\tag 3.4
$$
is $(k_1+k_2-2c-1)$--cartesian.  It follows by comparison
of~(3.3) and~(3.4) that the canonical map $PO \to \Phi(X_3)$ is
$(k_1+k_2-2c)$--connected (when $2c+1 \ge c$).  Applying $\Psi$ to~(3.4),
we obtain a $(k_1+k_2-3c)$--cartesian square.  The map $\Psi(PO)
\to \Psi\Phi(X_3)$ is $(k_1+k_2-3c)$--connected (when $2\kappa
\ge \kappa$), so the square obtained by applying $\Psi$ to~(3.3) is
$(k_1+k_2-3c-1)$--cartesian.
\qed
\enddemo

A weak homotopy functor $F \: \SS \to \SS$ satisfies the {\sl colimit
axiom\/} if it preserves filtered homotopy colimits up to weak homotopy.
This means that for any filtered diagram $X \: D \to \SS$ the canonical
map
$$
\hocolim_{d \in D} F(X_d) @>>> F(\hocolim_{d \in D} X_d)
\tag 3.5
$$
is a weak equivalence.  Any simplicial set is weakly equivalent to the
homotopy colimit of the filtered diagram of its finite sub-objects,
where a simplicial set is {\sl finite\/} if it has only finitely many
non-degenerate simplices.  Thus a functor satisfying the colimit axiom
is determined by its restriction to the subcategory of finite simplicial
sets.  Such functors are therefore also said to be {\sl finitary}.

Similarly, a functor $\Phi \: R(W, G) \to \SS$ satisfies the {\sl
colimit axiom\/} if it preserves filtered homotopy colimits up to weak
equivalence.  An object of $R(W, G)$ is said to be {\sl finite\/} if
it can be obtained from $W$ by attaching finitely many free $G$--cells.
Again a functor satisfying the colimit axiom is determined by its
restriction to the finite objects in $R(W, G)$.

The forgetful functor $u \: \SS/X \to \SS$ preserves filtered homotopy
colimits.  Hence if $F \: \SS \to \SS$ satisfies the colimit axiom,
then so does the composite functor $\Phi = F \circ u \: \SS/X \to \SS$.

\remark{Remark 3.6}
Let $\SS^\lambda\!/X$ be the full subcategory of $\SS/X$ with objects the
$\lambda$--connected maps $X' \to X$.  The conditions $E_1(c, \kappa)$ and
$E_2(c, \kappa)$ then make sense for functors $\Phi \: \SS^\lambda\!/X \to
\SS$, and all of the results of this paper also apply to functors with
such a restricted domain of definition.  One could even consider {\sl
germs} of functors $\SS/X \to \SS$, ie, equivalence classes of functors
$\Phi \: \SS^\lambda\!/X \to \SS$ defined for some integer $\lambda$,
with two such functors $\Phi$ and $\Phi' \: \SS^{\lambda'}\!/X \to \SS$
considered to be equivalent if there is a $\lambda''$ such that $\Phi |
\SS^{\lambda''}\!/X = \Phi' | \SS^{\lambda''}\!/X$.  For simplicity we
will not include this extra generality in our notation.
\endremark

\head Excisive approximation \endhead

If $X' \to X$ is an object of $\SS/X$, its {\sl fiberwise (unreduced)
cone} $C_X X'$ is the mapping cylinder $(X' \times \Delta^1) \cup_{X'} X$,
and its {\sl fiberwise (unreduced) suspension} $S_X X'$ is the union of
two such mapping cylinders along $X'$.  There is a cocartesian square of
simplicial sets over $X$:
$$
\xymatrix{
X' \ar[r] \ar[d] & C_X X' \ar[d] \\
C_X X' \ar[r] & S_X X'
}
$$
The functor $S_X$ increases the connectivity of simplicial sets and
maps by at least one.

Consider a weak homotopy functor $\Phi \: \SS/X \to \SS$.  Following
Goodwillie [\Cite{Go90}, \S1], we associate to $\Phi$ the weak homotopy
functor $T\Phi \: \SS/X \to \SS$ given by
$$
(T\Phi)(X') = \holim(\Phi(C_X X') \rightarrow \Phi(S_X X')
\leftarrow \Phi(C_X X')) \,.
$$
If $\Phi$ satisfies $E_1(c, \kappa)$ then $T\Phi$ satisfies $E_1(c,
\kappa-1)$.  There is a natural map $t\Phi \: \Phi \to T\Phi$.  Define
$T^n\Phi \: \SS/X \to \SS$ for $n\ge0$ by iteration, and let the weak
homotopy functor $P\Phi \: \SS/X \to \SS$ be the homotopy colimit $$
(P\Phi)(X') = \hocolim_n (T^n\Phi)(X') \,.
$$
Again there is a natural map $p\Phi \: \Phi \to P\Phi$, as functors
$\SS/X \to \SS$.  (Cf [\Cite{Go90}, 1.10].)  If $\Phi$ satisfies $E_1(c,
\kappa)$ then $P\Phi$ satisfies $E_1(c, \kappa')$ for all $\kappa'$.

We know that $P\Phi$ lifts to a functor $\SS/X \to \SS/P\Phi(X)$,
but we wish to modify it to a functor $\SS/X \to \SS/Y$, with
$Y = \Phi(X)$.  There is a commutative square
$$
\xymatrix{
\Phi(X') \ar[rr]^{p\Phi(X')} \ar[d] && P\Phi(X') \ar[d] \\
\Phi(X) \ar[rr]^{p\Phi(X)}_{\simeq} && P\Phi(X)
}
\tag 4.1
$$
induced by the unique morphism $X' \to X$ in $\SS/X$.  The lower
horizontal map is a weak equivalence by inspection of the construction
of $P\Phi$, using that $\Phi$ was assumed to preserve weak equivalences.
We set
$$
P'\Phi(X') = \holim(P\Phi(X') \rightarrow P\Phi(X) \leftarrow \Phi(X))
$$
equal to the homotopy limit (pullback) of the lower right hand part
of the diagram.  The commutative square~(4.1) then extends to
$$
\xymatrix{
\Phi(X') \ar[r]^-{p'\Phi(X')} \ar[d] & P'\Phi(X') \ar[r]^{\simeq} \ar[d]
& P\Phi(X') \ar[d] \\
\Phi(X) \ar[r]^{=} & \Phi(X) \ar[r]^{p\Phi}_{\simeq} & P\Phi(X) \,,
}
$$
where the right hand square only commutes up to homotopy.
Thus we can view $\Phi$ and $P'\Phi$ as functors $\SS/X \to \SS/Y$,
in which case there is a natural map $p'\Phi \: \Phi \to P'\Phi$.
Viewing $\Phi$ and $P'\Phi$ as functors to $\SS$, the natural map $p\Phi$
factors as
$$
p\Phi \: \Phi @>p'\Phi>> P'\Phi @>\simeq>> P\Phi \,,
\tag 4.2
$$
where the right hand map is a natural weak equivalence.  So $P'\Phi$ is a
weak homotopy functor, and if $\Phi$ satisfies $E_1(c, \kappa)$ for some
$\kappa$, then $P'\Phi$ satisfies $E_1(c, \kappa')$ for all $\kappa'$.

\remark{Remark 4.3}
Note that $P'\Phi(X)$ is typically not equal to $\Phi(X) = Y$, although
the canonical map $P'\Phi(X) \to Y$ is a weak equivalence, so $P'\Phi \:
\SS/X \to \SS/Y$ is not the canonical lift of its forgetful version $u
\circ P'\Phi \: \SS/X \to \SS$.
\endremark
\medskip

Suppose now that $\Phi \: \SS/X \to \SS$ satisfies condition $E_2(c,
\kappa)$, hence is stably excisive.  When the structure map $X' \to X$
is $k$--connected, for $k \ge \kappa$, it follows immediately that the
maps $t\Phi(X') \: \Phi(X') \to T\Phi(X')$, $p\Phi(X') \: \Phi(X')
\to P\Phi(X')$ and $p'\Phi(X') \: \Phi(X') \to P'\Phi(X')$ are all
$(2k-c)$--connected.

We say that two functors $\Phi, \Phi' \: \SS/X \to \SS$ satisfy condition
$\OO(c, \kappa)$ along a natural map $f \: \Phi \to \Phi'$ if whenever
$X' \to X$ is $k$--connected and $k \ge \kappa$, then $f(X') \: \Phi(X')
\to \Phi'(X')$ is $(2k-c)$--connected.  If $\Phi$ and $\Phi'$ satisfy
condition $\OO(c, \kappa)$ for some integers $c$ and $\kappa$, then we
say that $\Phi$ and $\Phi'$ {\sl agree to first order along $f$}.
(Cf [\Cite{Go90}, 1.13].)

\proclaim{Proposition 4.4}
Let $X$ be a simplicial set, $\Phi \: \SS/X \to \SS$ a stably excisive
weak homotopy functor, and let $Y = \Phi(X)$.  Then $P'\Phi \: \SS/X
\to \SS/Y$ is excisive, and $\Phi$ and $P'\Phi \: \SS/X \to \SS/Y$
agree to first order along $p'\Phi$.  If $\Phi$ is bounded below, then
so is $P'\Phi$.  If $\Phi$ satisfies the colimit axiom, then so
does $P'\Phi$.
\endproclaim

\demo{Proof}
Goodwillie proves in [\Cite{Go90}, 1.14] that the functor $P\Phi$ is
excisive and that $\Phi$ agrees with $P\Phi$ to first order along $p\Phi$.
In view of the weak equivalence in~(4.2), the same applies to $P'\Phi$.
We have noted above that if $\Phi$ is bounded below, then so are $T\Phi$,
$P\Phi$ and $P'\Phi$.  If $\Phi$ satisfies the colimit axiom, then so
does $T\Phi$, because filtered homotopy colimits commute with homotopy
pullbacks, up to weak equivalence.  (See theorem~1 on pages 215--216
in \cite{Ma71} for the corresponding statement for sets.)  Hence also
$P\Phi$ and $P'\Phi$ satisfy the colimit axiom, since the order of two
homotopy colimits can be commuted.
\qed
\enddemo

\head Goodwillie derivatives \endhead

Let $X$ be a simplicial set, let $\Phi \: \SS/X \to \SS$ be a weak
homotopy functor, and let $Y = \Phi(X)$.  As before, we may view $\Phi$
as a functor $\SS/X \to \SS/Y$ or $R(X) \to R(Y)$, without change in
notation.  Choose base points $x \in X$ and $y \in Y$.

There is a functor $i_0 = i_0(X, x) \: \SS_* \to R(X)$ that takes a based
simplicial set $T$ to the retractive simplicial set
$$
i_0(T) = X \vee_x T \,,
$$
where $r \: X \vee_x T \to X$ takes $T$ to the base point $x \in X$,
and $s \: X \to X \vee_x T$ is the standard inclusion.  This functor
preserves cocartesian squares.

There is a second functor $j_0 = j_0(Y, y) \: R(Y) \to \SS_*$ that takes a
retractive simplicial set $(Y', r, s)$ to the homotopy fiber
$$
j_0(Y') = \hofib_y(r \: Y' \to Y) \,,
$$
(with the natural base point that maps to $s(y) \in Y'$).  This functor
preserves $k$--cartesian squares for all $k$, hence also cartesian squares.

We shall later consider equivariant improvements $i$ and $j$ of $i_0$
and $j_0$, respectively, which may justify the notation.

If $\Phi$ is an excisive weak homotopy functor, then the composite functor
$$
L \: \SS_* @>i_0>> R(X) @>\Phi>> R(Y) @>j_0>> \SS_*
$$
is an excisive weak homotopy functor that takes $*$ to $L(*) =
\hofib_y(Y \to Y)$, which is contractible.  We say that $L$ is a
{\sl linear\/} functor.  It corresponds to a generalized (reduced)
homology theory given by $L_*(T) = \pi_*(L(T))$, with an associated
coefficient spectrum $\bold L = \{ n \mapsto L(S^n) \}$ (modulo a
technical rectification, as in [\Cite{Go90}, 0.1]).  There is a natural
weak equivalence $\Omega^\infty(\bold L \wedge T) \to L(T)$, at least for
finite simplicial sets $T$.  (Cf \cite{Go90} and proposition~10.4
below.)

Even if $\Phi$ is not excisive, we can still form the composite functor
$j_0 \circ \Phi \circ i_0$ and assemble the based simplicial sets $(j_0
\circ \Phi \circ i_0)(S^n)$ into a spectrum.  For any weak homotopy
functor $\Phi \: \SS/X \to \SS$ let
$$
\partial^x_y \Phi(X)_n = \hofib_y(\Phi(X \vee_x S^n) \to \Phi(X))
\tag 5.1
$$
for $n\ge0$.  There is a natural chain of maps
$$
\align
\partial^x_y \Phi(X)_n
&@<\simeq<<
\hofib_y(\Phi(X \vee_x S^n) \to \Phi(X \vee_x CS^n)) \\
&\longrightarrow
\hofib_y(\Phi(X \vee_x CS^n) \to \Phi(X \vee_x S^{n+1}))
@<\simeq<<
\Omega \partial^x_y \Phi(X)_{n+1}
\endalign
$$
where the second map is the natural one between the horizontal homotopy
fibers in the commutative diagram
$$
\xymatrix{
\Phi(X \vee_x S^n) \ar[r] \ar[d] &
\Phi(X \vee_x CS^n) \ar[d] \\
\Phi(X \vee_x CS^n) \ar[r] &
\Phi(X \vee_x S^{n+1})
}
$$
and the other two maps are natural weak equivalences derived from the
weak equivalence $X \vee_x CS^n \to X$ and the Puppe sequence.  We let
$$
\partial^x_y \Phi(X) = \{n \mapsto \partial^x_y \Phi(X)_n\}
\tag 5.2
$$
be the spectrum obtained from this sequence of based simplicial sets and
(weak) adjoint structure maps by the functorial rectification procedure
of [\Cite{Go90}, 0.1].  By definition, $\partial^x_y \Phi(X)$ is the {\sl
Goodwillie derivative\/} of $\Phi$ at $X$, with respect to the base
points $x \in X$ and $y \in Y = \Phi(X)$.  (Cf [\Cite{Go90}, 1.16].)

A natural map $f \: \Phi \to \Phi'$ of functors $\SS/X \to \SS$ induces
maps $\partial f \: \partial^x_y \Phi(X)_n \to \partial^x_{y'} \Phi'(X)_n$
for all $n\ge0$, and a spectrum map $\partial f \: \partial^x_y \Phi(X)
\to \partial^x_{y'} \Phi'(X)$.  This presupposes that $\Phi'(X)$ is
given the base point $y' = f(X)(y)$, where $y$ is the chosen base point
in $\Phi(X)$ and $f(X) \: \Phi(X) \to \Phi'(X)$.

\proclaim{Proposition 5.3}
If $\Phi$ and $\Phi' \: \SS/X \to \SS$ agree to first order along $f \:
\Phi \to \Phi'$, then $f$ induces a meta-stable equivalence of spectra
$\partial f \: \partial^x_y \Phi(X) \to \partial^x_{y'} \Phi'(X)$.
\endproclaim

\demo{Proof}
This is basically [\Cite{Go90}, 1.17].  Suppose that $\Phi$ and $\Phi'$
satisfy $\OO(c, \kappa)$.  The retraction $X \vee_x S^n \to X$ is
$n$--connected, so for $n \ge \kappa$ the map $f(X \vee_x S^n) \: \Phi(X
\vee_x S^n) \to \Phi'(X \vee_x S^n)$ is $(2n-c)$--connected.  In a similar way
$f(X) \: \Phi(X) \to \Phi'(X)$ is a weak equivalence.  Hence the map of
homotopy fibers $\partial f \: \partial^x_y \Phi(X)_n \to \partial^x_{y'}
\Phi'(X)_n$ is $(2n-c)$--connected.
\qed
\enddemo

Let $X$ be a simplicial set, $\Phi \: \SS/X \to \SS$ a weak homotopy
functor, $Y = \Phi(X)$, and choose base points $x \in X$ and $y \in Y$.
Give $P'\Phi(X)$ (defined in \S4) the base point $y' = p'\Phi(X)(y)$.

\proclaim{Corollary 5.4}
If $\Phi$ is stably excisive, then $p'\Phi$ induces a meta-stable
equivalence of spectra $\partial(p'\Phi) \: \partial^x_y \Phi(X) \to
\partial^x_{y'} (P'\Phi)(X)$.
\endproclaim

When $F \: \SS \to \SS$ is a weak homotopy functor, $X$ a simplicial
set, $Y = F(X)$, $x \in X$, $y \in Y$ and $\Phi = F \circ u$, we let
$\partial^x_y F(X)_n = \partial^x_y \Phi(X)_n = \hofib_y(F(X \vee_x S^n)
\to F(X))$ and $\partial^x_y F(X) = \partial^x_y \Phi(X) = \{n \mapsto
\partial^x_y F(X)_n\}$.

\head Composite functors \endhead

Let $X$ be a simplicial set, $\Phi \: \SS/X \to \SS$ a functor, $Y =
\Phi(X)$, $\Psi \: \SS/Y \to \SS$, and $Z = \Psi(Y)$.  Suppose that
$\Phi$ and $\Psi$ are weak homotopy functors.  Let $Y_1 = P'\Phi(X)$,
$Z_1 = \Psi(Y_1)$ and $Z_2 = P'\Psi(Y_1)$.  Choose base points $x \in X$,
$y \in Y$ and $z \in Z$, and let $y_1 = p'\Phi(X)(x) \in Y_1$, $z_1 =
\Psi(p'\Phi(X))(z) \in Z_1$ and $z_2 = p'\Psi(Y_1)(z_1) \in Z_2$.
$$
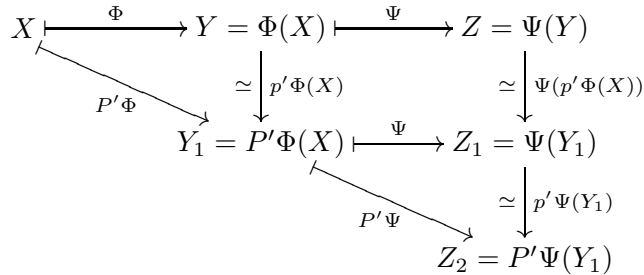
\xymatrix@!C{
{\phantom{XX}X} \ar@{|->}[r]^-{\Phi} \ar@{|->}[dr]_{P'\Phi} &
Y = \Phi(X) \ar@{|->}[r]^-{\Psi} \ar[d]^{p'\Phi(X)}_{\simeq} &
Z = \Psi(Y) \ar[d]^{\Psi(p'\Phi(X))}_{\simeq} \\
& Y_1 = P'\Phi(X) \ar@{|->}[r]^-{\Psi} \ar@{|->}[dr]_{P'\Psi} &
Z_1 = \Psi(Y_1) \ar[d]^{p'\Psi(Y_1)}_{\simeq} \\
&& Z_2 = P'\Psi(Y_1)
}
$$

\proclaim{Proposition 6.1}
Suppose that $\Phi$ and $\Psi$ are bounded below, stably excisive
functors.  Then the composite functors $\Psi \circ \Phi$ and $P'\Psi
\circ P'\Phi \: \SS/X \to \SS$ agree to first order along $p'\Psi
\circ \Psi(p'\Phi) = P'\Psi(p'\Phi) \circ p'\Psi$.
\endproclaim

\demo{Proof}
Assume that $\Phi$ and $\Psi$ satisfy $E_1(c, \kappa)$ and $E_2(c,
\kappa)$, for some sufficiently large integers $c$ and~$\kappa$.  Let $X'
\to X$ be a $k$--connected map, with $k \ge \kappa+c$.  Then $p'\Phi(X')
\: \Phi(X') \to P'\Phi(X')$ is $(2k-c)$--connected, so $\Psi(p'\Phi)(X') \:
(\Psi \circ \Phi)(X') \to (\Psi \circ P'\Phi)(X')$ is $(2k-2c)$--connected.
Furthermore, $P'\Phi(X') \to P'\Phi(X)$ is $(k-c)$--connected, as noted
after diagram~(4.2), so $p'\Psi(P'\Phi(X')) \: (\Psi \circ P'\Phi)(X')
\to (P'\Psi \circ P'\Phi)(X')$ is $(2k-3c)$--connected.  Thus $p'\Psi
\circ \Psi(p'\Phi)$ satisfies $\OO(3c, \kappa+c)$, and $\Psi \circ \Phi$
and $P'\Psi \circ P'\Phi \: \SS/X \to \SS$ agree to first order.
\qed
\enddemo

Recall also that for $\Phi$ and $\Psi$ bounded below and stably excisive
the composite functor $\Psi \circ \Phi \: \SS/X \to \SS$ is bounded
below and stably excisive (proposition~3.2), hence agrees to first
order with $P'(\Psi \circ \Phi) \: \SS/X \to \SS$ (proposition~4.4).
We are therefore legitimately interested in its derivative $\partial^x_z
(\Psi \circ \Phi)(X)$.

\proclaim{Proposition 6.2}
Suppose that $\Phi$ and $\Psi$ are bounded below, stably excisive
functors.  Then there are natural meta-stable
equivalences
\roster
\item
$\partial(p'\Phi) \: \partial^x_y \Phi(X) \simeq \partial^x_{y_1}
	P'\Phi(X)$,
\item
$\partial(p'\Psi) \: \partial^{y_1}_{z_1} \Psi(Y_1) \simeq
	\partial^{y_1}_{z_2} (P'\Psi)(Y_1)$,
\item
$\partial(p'\Psi \circ \Psi(p'\Phi)) \: \partial^x_z (\Psi \circ \Phi)(X)
	\simeq \partial^x_{z_2} (P'\Psi \circ P'\Phi)(X)$
\endroster
and a strict equivalence
\roster
\itemold{$(4)$}
$\partial^y_z \Psi(Y) \simeq \partial^{y_1}_{z_1} \Psi(Y_1)$.
\endroster
\endproclaim

\demo{Proof}
By propositions~4.4 and~6.1, the pairs of functors $\Phi$ and $P'\Phi$,
$\Psi$ and $P'\Psi$, and the functors $\Psi \circ \Phi$ and $P'\Psi \circ
P'\Phi$ agree to first order, respectively.  Hence their derivatives
are meta-stably equivalent by proposition~5.3.

Case~(4) remains.  There is a commutative square in $\SS/Y$
$$
\xymatrix{
Y \vee_y S^n \ar[d] \ar[r]^{\simeq} & Y_1 \vee_{y_1} S^n \ar[d] \\
Y \ar[r]^{p'\Phi(X)}_\simeq & Y_1
}
$$
where the vertical maps take $S^n$ to the respective base points, and the
upper horizontal map is the identity on $S^n$.  The lower horizontal
map is a weak equivalence, as in remark~4.3, hence so is the upper
horizontal map.  Applying $\Psi$ and taking vertical homotopy fibers
yields a weak equivalence of $n$-th spaces
$$
\partial^y_z \Psi(Y)_n \to \partial^{y_1}_{z_1} \Psi(Y_1)_n \,.
$$
Thus the associated spectra are strictly equivalent.
\qed
\enddemo

\remark{Remark 6.3}
It follows that for the purpose of expressing the derivative of $\Psi
\circ \Phi$ in terms of the derivatives of $\Phi$ and $\Psi$, we are
free to replace the bounded below, stably excisive functors $\Phi$
and $\Psi$ by their bounded below, excisive approximations $P'\Phi$
and $P'\Psi$, respectively.  If $\Psi$ satisfies the colimit axiom,
then so does its replacement.

Equivalently, we may assume that $\Phi$ and $\Psi$ are themselves bounded
below, excisive functors.  Furthermore, the derivatives only depend on the
associated functors $\Phi \: R(X) \to R(Y)$, $\Psi \: R(X) \to R(Y)$
with composite $\Psi \circ \Phi \: R(X) \to R(Z)$.
\endremark

\head Multiple connected components \endhead

We now reduce to the case when $X$, $Y$ and $Z$ are connected.

Let $\Phi \: R(X) \to R(Y)$ and $\Psi \: R(Y) \to R(Z)$ be weak
homotopy functors, with $Y = \Phi(X)$ and $Z = \Psi(Y)$.
Write $X = \coprod_{\alpha \in A} X_\alpha$, $Y = \coprod_{\beta \in B}
Y_\beta$ and $Z = \coprod_{\gamma \in C} Z_\gamma$, where each $X_\alpha$,
$Y_\beta$ and $Z_\gamma$ is connected.  So $A = \pi_0(X)$, $B =
\pi_0(Y)$ and $C = \pi_0(Z)$.  Choose base points $x_\alpha \in
X_\alpha$, $y_\beta \in Y_\beta$ and $z_\gamma \in Z_\gamma$ for all
$\alpha$, $\beta$ and $\gamma$.

Let $\In_\alpha \: R(X_\alpha) \to R(X)$ be given by pushout along
$X_\alpha \subset X$, so $\In_\alpha(X'_\alpha) = X \cup_{X_{\alpha}}
X'_\alpha$.  Similarly let $\pr_\beta \: R(Y) \to R(Y_\beta)$ be given
by pullback along $Y_\beta \subset Y$, so $\pr_\beta(Y') = Y_{\beta}
\times_Y Y'$.  Let $\Phi^\alpha_\beta = \pr_\beta \circ \Phi \circ
\In_\alpha \: R(X_\alpha) \to R(Y_\beta)$.  Clearly $\In_\alpha$
preserves $k$--connected maps and ($k$--)cocartesian squares, while
$\pr_\beta$ preserves $k$--connected maps and ($k$--)cartesian squares.
So if $\Phi$ is bounded below, excisive, stably excisive or satisfies
the colimit axiom, then the same applies to $\Phi^\alpha_\beta$.

\proclaim{Lemma 7.1}
There is a natural strict equivalence
$$
\partial^{x_\alpha}_{y_\beta} \Phi(X) \simeq
\partial^{x_\alpha}_{y_\beta} \Phi^\alpha_\beta(X_\alpha) \,.
$$
\endproclaim

\demo{Proof}
Let $X'_\alpha = X_\alpha \vee_{x_\alpha} S^n$, $X' = \In_\alpha(X'_\alpha)
\cong X \vee_{x_\alpha} S^n$, $Y' = \Phi(X')$ and $Y'_\beta =
\pr_\beta(Y')$, so that $Y'_\beta = \Phi^\alpha_\beta(X'_\alpha)$.
The pullback square
$$
\xymatrix{
Y'_\beta \ar[r] \ar[d] & Y' \ar[d] \\
Y_\beta \ar[r] & Y
}
$$
is cartesian, so there is a weak equivalence
$$
\partial^{x_\alpha}_{y_\beta} \Phi(X)_n \cong \hofib_{y_\beta}(Y'
\to Y) \simeq \hofib_{y_\beta}(Y'_\beta \to Y_\beta) =
\partial^{x_\alpha}_{y_\beta} \Phi^\alpha_\beta(X_\alpha)_n \,.
\eqno{\sq}
$$

Let $\Psi^\beta = \Psi \circ \In_\beta \: R(Y_\beta) \to R(Z)$ and
consider a base point $z \in Z$.  For each $Y'$ in $R(Y)$ and $\beta
\in B$ let $Y'_\beta = \pr_\beta(Y')$ in $R(Y_\beta)$.

\proclaim{Proposition 7.2}
Let $\Psi \: R(Y) \to R(Z)$ be a bounded below, excisive functor
that satisfies the colimit axiom.  Then the functors $Y' \mapsto
\hofib_z(\Psi(Y') \to Z)$ and
$$
Y' \mapsto \bigvee_{\beta \in B} \hofib_z(\Psi^\beta(Y'_\beta) \to Z)
$$
agree to first order along a natural chain of maps.
\endproclaim

\demo{Proof}
The retraction $Y' \to \In_\beta(Y'_\beta)$ induces a retraction $\Psi(Y')
\to \Psi^\beta(Y'_\beta)$.  More generally, for each finite subset $S
\subset B$ let
$$
Y^*_S = \coprod_{\beta\in S} Y'_\beta \cup
\coprod_{\beta\notin S} Y_\beta \,.
$$
There is then a strongly cocartesian (cf [\Cite{Go92}, 2.1]) cubical
diagram
$$
(T \subset S) \longmapsto Y^*_{S \setminus T}
$$
in $R(Y)$.  Applying the excisive functor $\Psi$ yields a strongly
cartesian cubical diagram, where each map admits a section.
Hence there is a weak equivalence
$$
\hofib_z(\Psi(Y^*_S) \to Z)
@>\simeq>>
\prod_{\beta \in S} \hofib_z(\Psi^\beta(Y'_\beta) \to Z) \,.
$$
Passing to homotopy colimits over $S \subset B$, and using that
$\Psi$ satisfies the colimit axiom, yields a weak equivalence
$$
\hofib_z(\Psi(Y') \to Z)
@>\simeq>>
\hocolim_{S \subset B} \prod_{\beta \in S} \hofib_z(\Psi^\beta(Y'_\beta)
\to Z) \,.
$$

When $\Psi$ satisfies $E_1(c, \kappa)$ and $Y' \to Y$ is $k$--connected,
with $k \ge \kappa$, then $\In_\beta(Y'_\beta) \to Y$ is $k$--connected and
$\Psi^\beta(Y'_\beta) \to Z$ is $(k-c)$--connected, for each $\beta$.  So
each space $\hofib_z(\Psi^\beta(Y'_\beta) \to Z)$ is $(k-c-1)$--connected,
each inclusion
$$
\bigvee_{\beta \in S} \hofib_z(\Psi^\beta(Y'_\beta) \to Z)
@>>>
\prod_{\beta \in S} \hofib_z(\Psi^\beta(Y'_\beta) \to Z)
$$
is $(2k-2c-1)$--connected, and the resulting map
$$
\hocolim_{S \subset B} \bigvee_{\beta \in S} \hofib_z(\Psi^\beta(Y'_\beta)
\to Z)
@>>>
\hocolim_{S \subset B} \prod_{\beta \in S} \hofib_z(\Psi^\beta(Y'_\beta)
\to Z)
$$
is $(2k-2c-1)$--connected.  The source of this map is naturally equivalent
to $\bigvee_{\beta \in B} \hofib_z(\Psi^\beta(Y'_\beta) \to Z)$.
\qed
\enddemo

Let $\Phi^\alpha = \Phi \circ \In_\alpha \: R(X_\alpha) \to R(Y)$
and $\Psi_\gamma = \pr_\gamma \circ \Psi \: R(Y) \to R(Z_\gamma)$.

\proclaim{Proposition 7.3}
Let $\Phi \: R(X) \to R(Y)$ and $\Psi \: R(Y) \to R(Z)$ be bounded below,
and suppose that $\Psi$ is excisive and satisfies the colimit axiom.
There is a natural chain of meta-stable equivalences
$$
\partial^{x_\alpha}_{z_\gamma}(\Psi_\gamma \circ \Phi^\alpha)(X_\alpha)
\simeq
\bigvee_{\beta \in B}
\partial^{x_\alpha}_{z_\gamma}(\Psi^\beta_\gamma \circ
\Phi^\alpha_\beta)(X_\alpha) \,.
$$
\endproclaim

\demo{Proof}
We keep the notation introduced in the last two proofs.  Then
$$
\align
\partial^{x_\alpha}_{z_\gamma}(\Psi_\gamma \circ \Phi^\alpha)(X_\alpha)_n
&= \hofib_{z_\gamma}((\Psi_\gamma \circ \Phi^\alpha)(X'_\alpha)
\to Z_\gamma) \\
&= \hofib_{z_\gamma}(\Psi_\gamma(Y') \to Z_\gamma) \\
&\sim \bigvee_{\beta \in B}
\hofib_{z_\gamma}(\Psi^\beta_\gamma(Y'_\beta) \to Z_\gamma) \\
&= \bigvee_{\beta \in B}
\hofib_{z_\gamma}((\Psi^\beta_\gamma \circ \Phi^\alpha_\beta)(X'_\alpha)
\to Z_\gamma) \\
&= \bigvee_{\beta \in B}
\partial^{x_\alpha}_{z_\gamma}(\Psi^\beta_\gamma \circ \Phi^\alpha_\beta)
(X_\alpha)_n
\endalign
$$
where $\sim$ denotes agreement up to first order as functors of $Y'$,
by proposition~7.2.  Since $\Phi$ is bounded below, the map $Y' \to Y$
is $(n-c)$--connected for some constant~$c$, and so this chain of maps is
$(2n-c)$--connected for some (other) constant~$c$.  Hence the associated
spectrum map is a meta-stable equiva-lence.\qed
\enddemo

\head The Kan loop group \endhead

By the results of section~6 we need only consider bounded below, excisive
functors $\Phi \: R(X) \to R(Y)$ and $\Psi \: R(Y) \to R(Z)$, where $X$
is a simplicial set, $Y = \Phi(X)$ and $Z = \Psi(Y)$, and by section~7 we
may assume that $X$, $Y$ and $Z$ are all connected.  Choose base points
$x \in X$, $y \in Y$ and $z \in Z$, and recall the functors $i_0 = i_0(X,
x)$ and $j_0 = j_0(Z, z)$ from \S5.  We are interested in the coefficient
spectrum of the composite functor:
$$
\SS_* @>i_0>> R(X) @>\Phi>> R(Y) @>\Psi>> R(Z) @>j_0>> \SS_* \,.
$$
We would like to be able to express this in terms of the coefficient
spectra of the composites $j_0 \Phi i_0 \: \SS_* \to R(X) \to R(Y) \to
\SS_*$ and $j_0 \Psi i_0 \: \SS_* \to R(Y) \to R(Z) \to \SS_*$.  But,
as consideration of the two special cases $\Phi = P_X(\id)$ and $\Psi =
P_Y(\id)$ indicates, where $\id \: \SS \to \SS$ is the identity functor,
this is not likely to be possible.  Composition with the functors $j_0
\: R(Y) \to \SS_*$ and $i_0 \: \SS_* \to R(Y)$ does not retain enough
information.

We shall instead replace the category $\SS_* = R(*)$ by the category
$R(*, \Omega_y(Y))$, where $\Omega_y(Y)$ is the Kan loop group of $(Y,
y)$, and replace the functors $i_0$ and $j_0$ by suitable inverse weak
equivalences $i \: R(*, \Omega_y(Y)) \to R(Y)$ and $j \: R(Y) \to R(*,
\Omega_y(Y))$, respectively.  The Goodwillie derivative $\partial^x_y
\Phi(X)$ then becomes a free left $\Omega_y(Y)$--spectrum.  Furthermore,
there is a model for the Goodwillie derivative $\partial^y_z \Psi(Y)$
that becomes a spectrum with right $\Omega_y(Y)$--action.  It turns
out that these extra simplicial group actions suffice to express
the derivative $\partial^x_z (\Psi \circ \Phi)(X)$ in terms of these
equivariant Goodwillie derivatives, leading to the chain rule.

\medskip

Suppose now that $Y$ is a connected simplicial set, with a chosen
base point $y \in Y$.  The {\sl Kan loop group\/} of $Y$ is then a
functorially defined simplicial group $\Omega_y(Y)$.  (Cf \cite{Wa96},
where the Kan loop group $\Omega_y(Y)$ is denoted $G(Y)$.)  There is a
principal $\Omega_y(Y)$--bundle
$$
\Omega_y(Y) @>>> \tilde\Gamma(Y) @>\pi>> \Gamma(Y)
$$
with $\tilde\Gamma(Y)$ weakly contractible, and a natural inclusion
$\gamma \: Y \to \Gamma(Y)$ which is a weak equivalence.  There are
natural base points in $\Gamma(Y)$ and $\tilde\Gamma(Y)$, and $\pi$
and $\gamma$ preserve these points.  Let $\epsilon \: \tilde\Gamma(Y)
\to *$ be the (unique) weak equivalence.

\proclaim{Proposition 8.1}
There are natural pointed weak homotopy functors
$$
\align
i &\: R(*, \Omega_y(Y)) \to R(Y) \\
j &\: R(Y) \to R(*, \Omega_y(Y)) \,,
\endalign
$$
such that $i = i(Y, y)$ preserves $k$--connected maps and ($k$--)cocartesian
squares, and $j = j(Y, y)$ preserves $k$--connected maps and
($k$--)cartesian squares.  There is a natural weak equivalence from the
identity functor on $R(Y)$ to the composite $i \circ j \: R(Y) \to R(Y)$.
\endproclaim

\remark{Remark 8.2}
This is closely related to [\Cite{Wa85}, 2.1.4].  There is also a natural
chain of weak equivalences from the composite $j \circ i \: R(*,
\Omega_y(Y)) \to R(*, \Omega_y(Y))$ to the identity functor on $R(*,
\Omega_y(Y))$, but we shall not make any use of it in this paper.
\endremark

\demo{Proof}
We construct $i$ as a composite
$$
i \: R(*, \Omega_y(Y)) @>\epsilon^*>> R(\tilde\Gamma(Y), \Omega_y(Y))
@>\pi_!>> R(\Gamma(Y)) @>R\gamma^*>> R(Y)
$$
and $j$ as a composite
$$
j \: R(Y) @>\gamma_*>> R(\Gamma(Y)) @>\pi^!>> R(\tilde\Gamma(Y),
\Omega_y(Y)) @>\epsilon_*>> R(*, \Omega_y(Y)) \,.
$$

The pointed functors $\epsilon_*$ and $\epsilon^*$ are given by
pushout and pullback along the map $\epsilon \: \tilde\Gamma(Y) \to *$,
respectively.  These are weak homotopy functors because the structural
section $s$ is always a cofibration, and the map $\epsilon$ is a
weak equivalence.  On an object $(E, r, s)$ of $R(\tilde\Gamma(Y),
\Omega_y(Y))$ the composite $\epsilon^* \epsilon_*(E)$ equals
$\tilde\Gamma(Y) \times (* \cup_{\tilde\Gamma(Y)} E)$, and there is
a natural weak equivalence from the identity on $R(\tilde\Gamma(Y),
\Omega_y(Y))$ to $\epsilon^* \circ \epsilon_*$ in view of the commutative
diagram:
$$
\xymatrix{
E \ar[rr]^-{\simeq} \ar[dr]_{\simeq} && {\tilde\Gamma(Y)} \times
(* \cup_{\tilde\Gamma(Y)} E) \ar[dl]^{\simeq}_{\pr_2} \\
& {*} \cup_{\tilde\Gamma(Y)} E 
}
$$
Here the lower left arrow collapses the image of $s \: \tilde\Gamma(Y)
\to E$ to a point, while the upper horizontal arrow has components $r \:
E \to \tilde\Gamma(Y)$ and the map just mentioned.

The pointed functors $\pi_!$ and $\pi^!$ are given by passage to
$\Omega_y(Y)$--orbits and pullback along $\pi \: \tilde\Gamma(Y)
\to \Gamma(Y)$, respectively.  These are weak homotopy functors
because objects of $R(\tilde\Gamma(Y), \Omega_y(Y))$ are free
$\Omega_y(Y)$--simplicial sets, and the bundle projection $\pi$ is
a Kan fibration.  Let $(W, r, s)$ be an object of $R(\Gamma(Y))$.
The composite value $\pi_! \pi^!(W) = (\tilde\Gamma(Y) \times_{\Gamma(Y)}
W)/\Omega_y(Y)$ is naturally isomorphic to $W$, so there is a natural
isomorphism from the identity on $R(\Gamma(Y))$ to $\pi_! \circ \pi^!$.

The pointed functor $\gamma_*$ is given by pushout along $\gamma \:
Y \to \Gamma(Y)$.  It is a weak homotopy functor because $\gamma$ is
a cofibration.  The construction of the pointed functor $R\gamma^*$
is a little more complicated, and we will do it in two steps.  Let $(W,
r, s)$ be an object of $R(\Gamma(Y))$.  We first define $\gamma^\#(W)$
in $\SS/Y$ as the homotopy pullback
$$
\xymatrix{
\gamma^\#(W) \ar[r]^\simeq \ar[d] & W \ar[d]^r \\
Y \ar@{ >->}[r]^\gamma_\simeq & \Gamma(Y)
}
$$
along $\gamma$.  This defines a weak homotopy functor $\gamma^\# \:
R(\Gamma(Y)) \to \SS/Y$, because we take homotopy pullback rather than
pullback.  By functoriality, $\gamma^\#(W)$ contains
$\gamma^\#(\Gamma(Y))$ as a retract in $\SS/Y$.  We next define
$R\gamma^*(W)$ in $R(Y)$ as the pushout
$$
\xymatrix{
\gamma^\#(\Gamma(Y)) \ar@{ >->}[rr]^{\gamma^\#(s)} \ar[d]^{\simeq} &&
\gamma^\#(W) \ar[d]^{\simeq} \\
Y \ar@{ >->}[rr] && R\gamma^*(W)
}
$$
in $\SS/Y$.  This defines a weak homotopy functor $R\gamma^* \:
R(\Gamma(Y)) \to R(Y)$, because the section $\gamma^\#(s)$ is a
(split) cofibration.  It is pointed by inspection.

We claim that there is a natural weak equivalence from the identity on
$R(Y)$ to $R\gamma^* \circ \gamma_*$.  Let $(Y', r, s)$ be an object of
$R(Y)$ and write $W = \gamma_*(Y') = \Gamma(Y) \cup_Y Y'$.  The canonical
inclusion $Y' \to W$ and the retraction $r \: Y' \to Y$ taken together
define a map from $Y'$ to the (strict) pullback in the diagram defining
$\gamma^\#(W)$.  Continuing by the canonical map from the pullback
to the homotopy pullback defines a natural map $Y' \to \gamma^\#(W)$.
It is a weak equivalence, in view of the commutative diagram:
$$
\xymatrix{
Y' \ar[rr]^\simeq \ar@{ >->}[dr]_\simeq &&
\gamma^\#(W) \ar[dl]^\simeq \\
& W 
}
$$
A diagram chase shows that the composite weak equivalence $Y' \to
\gamma^\#(W) \to R\gamma^*(W)$ is a morphism in $R(Y)$.  Hence this
defines the desired natural weak equivalence.

It is clear by inspection that $\epsilon^*$, $\pi_!$, $\gamma^\#$
and $R\gamma^*$ preserve $k$--connected maps and cartesian squares, and
that $\epsilon_*$, $\pi^!$ and $\gamma_*$ preserve $k$--connected maps
and cocartesian
squares.
\qed
\enddemo

\proclaim{Proposition 8.3}
The functor $i_0 \: \SS_* \to R(Y)$ factors up to a natural chain of
weak equivalences as the composite
$$
\SS_* @>\Omega_y(Y)_+ \wedge (-)>> R(*, \Omega_y(Y)) @>i>> R(Y) \,,
$$
where the left hand functor takes $T$ to $\Omega_y(Y)_+ \wedge T$.

The functor $j_0 \: R(Y) \to \SS_*$ factors up to a natural chain of
weak equivalences as the composite
$$
R(Y) @>j>> R(*, \Omega_y(Y)) @>w>> \SS_* \,,
$$
where the right hand functor forgets the $\Omega_y(Y)$--action.
\endproclaim

\demo{Proof}
Let $T$ be a based simplicial set.  We interpret
$\Omega_y(Y)_+ \wedge T$ as the pushout of
$$
* \leftarrow \Omega_y(Y) \times * \rightarrow \Omega_y(Y) \times T
$$
and apply $\pi_! \epsilon^*$ to obtain the pushout square:
$$
\xymatrix{
{\tilde\Gamma(Y)} \ar[r] \ar[d]^{\pi} &
{\tilde\Gamma(Y)} \times T \ar[d] \\
\Gamma(Y) \ar[r] &
(\tilde\Gamma(Y) \times (\Omega_y(Y)_+ \wedge T))/\Omega_y(Y)
}
$$
Using the weak equivalence $\gamma \: Y \to \Gamma(Y)$, and the inclusion
of the canonical base point into $\tilde\Gamma(Y)$, we obtain a natural
map from the pushout square
$$
\xymatrix{
{*} \ar[r] \ar[d] & T \ar[d] \\
Y \ar[r] & i_0(T)
}
$$
to the square above.  It is clearly a weak equivalence at the three
upper or left hand corners, hence the pushout map
$$
i_0(T) = Y \vee_y T @>\simeq>> \pi_! \epsilon^*(\Omega_y(Y)_+ \wedge T) \,.
$$
is also a natural weak equivalence.  The natural weak equivalences
$\gamma^\#(W) \to W$ and $\gamma^\#(W) \to R\gamma^*(W)$ from the two
diagrams defining $R\gamma^*$, and the factorization $i = R\gamma^*
\circ \pi_! \epsilon^*$, provide the remaining chain of natural weak
equivalences linking $i_0(T)$ to $i(\Omega_y(Y)_+ \wedge T)$.

Next, let $(Y', r, s)$ be a retractive simplicial set in $R(Y)$.
We evaluate $j = \epsilon_* \pi^! \gamma_*$ on $Y'$ and obtain a commutative
diagram in $\SS_*$:
$$
\xymatrix{
Y' \ar[d] \ar[r]^-{\simeq} & \Gamma(Y) \cup_Y Y' \ar[d] &
{\tilde\Gamma(Y)} \times_{\Gamma(Y)} (\Gamma(Y) \cup_Y Y') \ar@{->>}[l]
\ar[d] \ar[r]^-{\simeq} & j(Y') \ar[d] \\
Y \ar[r]^{\simeq}_{\gamma} & \Gamma(Y) &
{\tilde\Gamma(Y)} \ar@{->>}[l]^{\pi} \ar[r]^\simeq_{\epsilon}
& {*}
}
$$
The left and right hand squares have horizontal weak equivalences, while
the middle square is cartesian, using again that $\pi$ is a Kan fibration.
Hence the induced maps of vertical homotopy fibers at the canonical base
points of $Y$, $\Gamma(Y)$, $\tilde\Gamma(Y)$ and $*$ define a natural
chain of weak equivalences linking $j_0(Y') = \hofib_y(Y' \to Y)$ to
$\hofib(j(Y') \to *)$, and thus to $j(Y')$, as based simplicial sets.
\qed
\enddemo

\proclaim{Lemma 8.4}
Let $\Phi \: R(X) \to R(Y)$ be a weak homotopy functor, and let $\phi
= j \circ \Phi \circ i \: R(*, \Omega_x(X)) \to R(*, \Omega_y(Y))$.
If $\Phi$ satisfies the colimit axiom, then so does $\phi$.
\endproclaim

\demo{Proof}
By a straightforward inspection of the definitions, the functors $i
= R\gamma^* \circ \pi_! \circ \epsilon^*$ and $j = \epsilon_* \circ
\pi^! \circ \gamma_*$ preserve filtered homotopy colimits up to weak
equivalence.
\qed
\enddemo

\head Equivariant derivatives \endhead

We now recast the Goodwillie derivative in an equivariant setting.

Let $X$ be a simplicial set, let $\Phi \: \SS/X \to \SS$ be a weak
homotopy functor, and let $Y = \Phi(X)$.  Consider $\Phi$ as a pointed
functor $\Phi \: R(X) \to R(Y)$.  We choose base points $x \in X$ and $y
\in Y$, and shall first assume that $X$ and $Y$ are connected.  Let $H =
\Omega_x(X)$ and $G = \Omega_y(Y)$ be the Kan loop groups of $(X, x)$
and $(Y, y)$, respectively.

Recall the functors $i = i(X, x) \: R(*, H) \to R(X)$ and $j = j(Y,
y) \: R(Y) \to R(*, G)$.  Let $\phi \: R(*, H) \to R(*, G)$ be the
composite functor
$$
\phi \: R(*, H) @>i>> R(X) @>\Phi>> R(Y) @>j>> R(*, G)
\,.
\tag 9.1
$$
By proposition~8.1, $\phi$ is also a pointed weak homotopy functor.
(Here pointed means that $\phi(*) = *$.)

The categories $R(*, H)$ and $R(*, G)$ are in fact {\sl enriched\/} over
the category $\SS_*$ of based simplicial sets.  Given objects $U$ and
$V$ in $R(*, H)$, the (based) {\sl simplicial mapping space} $\Map_H(U,
V)$ is the based simplicial set with $p$--simplices the set of morphisms
$\Delta^p_+ \wedge U \to V$ in $R(*, H)$.  The usual set of morphisms $U
\to V$ in $R(*, H)$ can be recovered as the $0$--simplices of $\Map_H(U,
V)$.  Similarly, $\Map_G(-, -)$ is the simplicial mapping space in the
category $R(*, G)$.  A pointed {\sl simplicial functor} $\phi \: R(*,
H) \to R(*, G)$ comes equipped with a based map
$$
\Map_H(U, V) @>\phi>> \Map_G(\phi(U), \phi(V))
$$
which on $0$--simplices takes a map $f \: U \to V$ to the usual image
$\phi(f) \: \phi(U) \to \phi(V)$.

Any (pointed) weak homotopy functor $\phi$ can be promoted to a weakly
equivalent (pointed) simplicial weak homotopy functor $\check\phi$,
following [\Cite{Wa85}, 3.1].  For a based, free $H$--simplicial set $W$, let
$W^{\Delta^q} = \Map_H(H_+ \wedge \Delta^q_+, W)$.  Then $W^{\Delta^q}$
is again a based, free $H$--simplicial set.  The functor $[q] \mapsto
\phi(W^{\Delta^q})$ is a simplicial object in based, free $G$--simplicial
sets.  Let
$$
\check\phi(W) = \diag( [q] \mapsto \phi(W^{\Delta^q}) )
= \coprod_{q\ge0} (\phi(W^{\Delta^q}) \times \Delta^q) / \sim
$$
be the associated diagonal based, free $G$--simplicial set.  Then
$\check\phi$ is naturally a simplicial functor.  The required map
$$
\Map_H(U, V) @>\check\phi>> \Map_G(\check\phi(U), \check\phi(V))
$$
takes each $p$--simplex $f \: \Delta^p_+ \wedge U \to V$ to a $p$--simplex
$\check\phi(f) \: \Delta^p_+ \wedge \check\phi(U) \to \check\phi(V)$.
The latter is the simplicial map given in degree $q$ by the map
$$
\check\phi(f)_q \:
(\Delta^p_q)_+ \wedge \phi(U^{\Delta^q})
\cong
\bigvee_{\alpha \: \Delta^q \to \Delta^p} \phi(U^{\Delta^q})
@>>>
\phi(V^{\Delta^q})
$$
that on the wedge summand indexed by $\alpha \: \Delta^q \to \Delta^p$
is the value of $\phi$ applied to the composite map
$$
U^{\Delta^q} @>(\alpha, \id)>> (\Delta^p_+ \wedge U)^{\Delta^q}
@>f^{\Delta^q}>> V^{\Delta^q} \,.
$$
Each projection $\Delta^q \to *$ induces a weak equivalence $W \to
W^{\Delta^q}$, and thus each $q$--fold degeneracy map $\phi(W) \to
\phi(W^{\Delta^q})$ is a weak equivalence.  Thus the inclusion of
$0$--simplices is a natural weak equivalence
$$
\phi @>\simeq>> \check\phi
\tag 9.2
$$
of weak homotopy functors $R(*, H) \to R(*, G)$, by the realization lemma.
We can therefore replace $\phi = j \circ \Phi \circ i$ by $\check\phi$
without changing its weak homotopy type.

To each based, free $H$--simplicial set $W$ we associate its {\sl cone} $CW
= W \wedge \Delta^1$, and its {\sl suspension} $\Sigma W = CW \cup_W CW =
W \wedge S^1$, where $S^1 = \Delta^1 \cup_{\partial \Delta^1} \Delta^1$.
By iteration, $\Sigma^n W = W \wedge S^n$, where $S^n = S^1 \wedge \dots
\wedge S^1$ ($n$ copies of $S^1$).

Using that $\check\phi$ is a pointed simplicial functor, we obtain a
natural based map
$$
\sigma \: \Sigma\check\phi(W) \to \check\phi(\Sigma W)
\tag 9.3
$$
in $R(*, G)$, as follows.  The identity map $W \wedge S^1 = \Sigma W$
is left adjoint to a based map
$$
S^1 \to \Map_H(W, \Sigma W) \,.
$$
Since $\check\phi$ is pointed simplicial there is a natural based map
$$
\Map_H(W, \Sigma W) @>\check\phi>> \Map_G(\check\phi(W), \check\phi(\Sigma W)) \,.
$$
The composite of these two maps is then right adjoint to the
desired map $\sigma$.

We consider the sequence of based, free $H$--simplicial sets $\Sigma^n H_+
= H_+ \wedge S^n$ in $R(*, H)$.  Applying $\check\phi$ we obtain a sequence
of based, free $G$--simplicial sets
$$
\partial \phi_n = \check\phi(\Sigma^n H_+) \,.
\tag 9.4
$$
The natural map $\sigma$ in the case of $W = \Sigma^n H_+$ then defines
the structure map from $\Sigma \partial \phi_n$ to $\partial
\phi_{n+1}$ in the free $G$--spectrum
$$
\partial \phi = \{n \mapsto \partial \phi_n\} \,.
\tag 9.5
$$
This is the {\sl equivariant Goodwillie derivative\/} of $\phi$.
By proposition~8.3 and the weak equivalence~(9.2), there is a natural
chain of weak equivalences linking the underlying based simplicial set of
$\partial \phi_n$ to $\partial^x_y \Phi(X)_n$.  Similarly the underlying
non-equivariant spectrum of $\partial \phi$ is strictly equivalent to
the (non-equivariant) Goodwillie derivative $\partial^x_y \Phi(X)$.
Hence $\partial \phi$ provides a model for the Goodwillie derivative as
a free left $G$--spectrum.

The simplicial enrichment ensures that the Goodwillie derivative $\partial
\phi$ also admits another simplicial group action, this time by $H$ acting
from the right.  The multiplication on $H$ defines a map $(\Sigma^n
H_+) \wedge H_+ \to (\Sigma^n H_+)$, which is left adjoint to a based map
$$
H_+ @>>> \Map_{H}(\Sigma^n H_+, \Sigma^n H_+) \,.
$$
Using that $\check\phi \: R(*, H) \to R(*, G)$ is a pointed simplicial
functor, we obtain a map of based simplicial sets
$$
\Map_H(\Sigma^n H_+, \Sigma^n H_+)
@>\check\phi>>
\Map_G(\check\phi(\Sigma^n H_+), \check\phi(\Sigma^n H_+)) \,.
$$
The composite of these two based maps is then right adjoint to a map
$$
\check\phi(\Sigma^n H_+) \wedge H_+ @>>> \check\phi(\Sigma^n H_+) \,,
$$
that defines the desired right action of $H$ on $\check\phi(\Sigma^n H_+)
= \partial \phi_n$ in the category of based, free $G$--simplicial sets.

The structure maps of $\partial \phi$ are likewise natural with respect
to this right $H$--action, hence $\partial \phi$ is also a spectrum with
right $H$--action.

Turning to the general case, when $X$ and $Y$ are not necessarily
connected, let $X_\alpha$ be the connected component of the base point $x$
in $X$ and let $Y_\beta$ be the connected component of the base point $y$
in $Y$.  Define $\Phi^\alpha_\beta \: R(X_\alpha) \to R(Y_\beta)$ as
in section~7, and let $\phi = j \circ \Phi^\alpha_\beta \circ i \: R(*,
\Omega_x(X_\alpha)) \to R(*, \Omega_y(Y_\beta))$.  Then we have just
seen that $\partial\phi \simeq \partial^x_y \Phi^\alpha_\beta(X_\alpha)$,
while $\partial^x_y \Phi^\alpha_\beta(X_\alpha) \simeq \partial^x_y
\Phi(X)$ by lemma~7.1.

We summarize this discussion in:

\definition{Definition 9.6}
Let $X$ be a simplicial set, $\Phi \: \SS/X \to \SS$ a weak homotopy
functor and $Y = \Phi(X)$.  Let $(X_\alpha, x)$ and $(Y_\beta, y)$
be based, connected components of $X$ and $Y$, with Kan loop groups
$\Omega_x(X_\alpha)$ and $\Omega_y(Y_\beta)$, respectively.  Let
$$
\phi = j \circ \Phi^\alpha_\beta \circ i \: R(*, \Omega_x(X_\alpha))
\to R(*, \Omega_y(Y_\beta))
\,,
$$
with $i = i(X_\alpha, x)$ and $j = j(Y_\beta, y)$, and let
$$
\check\phi(W) = | [q] \mapsto \phi(W^{\Delta^q}) |
$$
be the simplicial enrichment of $\phi$.  The {\sl equivariant Goodwillie
derivative\/} of $\phi$ is the free left $\Omega_y(Y_\beta)$--spectrum
with right $\Omega_x(X_\alpha)$--action:
$$
\partial \phi = \{n \mapsto \partial\phi_n = \check\phi(\Sigma^n
\Omega_x(X_\alpha)_+)\}
\,.
$$
Its underlying spectrum is strictly equivalent to the
(non-equivariant) Goodwillie derivative $\partial^x_y \Phi(X)$.
\enddefinition

\head Equivariant Brown--Whitehead representability\endhead

Let $H$ and $G$ be simplicial groups.  We say that a functor $\phi \:
R(*, H) \to R(*, G)$ is {\sl linear\/} if it is a pointed, excisive,
weak homotopy functor.  We show in this section that linear functors that
satisfy the colimit axiom are classified by their equivariant Goodwillie
derivative.  Recall from \S3 that an object $W$ of $R(*, H)$ is finite
if it can be obtained from $*$ by attaching finitely many free $H$--cells.

\proclaim{Proposition 10.1}
Let $\phi \: R(*, H) \to R(*, G)$ be a bounded below, linear functor.
There is a free left $G$--spectrum
$$
\partial \phi = \{n \mapsto \partial\phi_n = \check\phi(\Sigma^n H_+)\}
$$
with right $H$--action, and a natural map
$$
\theta(W) \: \partial \phi \wedge_H W @>>> \{n \mapsto \check\phi(\Sigma^n
W)\}
$$
of free $G$--spectra, which is a meta-stable equivalence for all finite
$W$ in $R(*, H)$.  If $\phi$ satisfies the colimit axiom, then the map
is a stable equivalence for all $W$ in $R(*, H)$.
\endproclaim

\demo{Proof}
The group action map $H_+ \wedge W \to W$ suspends to a map $\Sigma^n H_+
\wedge W \to \Sigma^n W$, which is left adjoint to a based map
$$
W @>>> \Map_H(\Sigma^n H_+, \Sigma^n W)
$$
of $H$--simplicial sets, where $H$ acts on the simplicial mapping space
by right multiplication in the domain.  Since $\check\phi$ is pointed
simplicial, there is a based map
$$
\Map_H(\Sigma^n H_+, \Sigma^n W) @>\check\phi>>
\Map_G(\check\phi(\Sigma^n H_+), \check\phi(\Sigma^n W))
$$
of $H$--simplicial sets.  The composite map is right adjoint to a map
$$
\check\phi(\Sigma^n H_+) \wedge_H W @>>> \check\phi(\Sigma^n W)
\tag 10.2
$$
of based, free $G$--simplicial sets.  These maps are compatible with the
spectrum structure maps for varying $n$, hence define the natural map
$\theta(W)$ of free $G$--spectra.

To finish the proof, we will use the following lemma:

\proclaim{Lemma 10.3}
Let
$$
\xymatrix{
W_0 \ar[r] \ar[d] & W_1 \ar[d] \\
W_2 \ar[r] & W_3
}
$$
be a cocartesian square in $R(*, H)$.  If $\theta(W_0)$, $\theta(W_1)$
and $\theta(W_2)$ are meta-stable equivalences, then so is $\theta(W_3)$.
\endproclaim

\demo{Proof}
Applying~(10.2) to the given cocartesian square yields a map from the square
$$
\xymatrix{
{\check\phi(\Sigma^n H_+)} \wedge_H W_0 \ar[r] \ar[d] &
{\check\phi(\Sigma^n H_+)} \wedge_H W_1 \ar[d] \\
{\check\phi(\Sigma^n H_+)} \wedge_H W_2 \ar[r] &
{\check\phi(\Sigma^n H_+)} \wedge_H W_3
}
$$
to the square:
$$
\xymatrix{
{\check\phi(\Sigma^n W_0)} \ar[r] \ar[d] &
{\check\phi(\Sigma^n W_1)} \ar[d] \\
{\check\phi(\Sigma^n W_2)} \ar[r] &
{\check\phi(\Sigma^n W_3)} 
}
$$
The first square is cocartesian with $(n-c)$--connected maps for some~$c$,
since $\Sigma^n H_+$ is $(n-1)$--connected, $\phi \simeq \check\phi$
is bounded below and each $W_i$ is $H$--free.  By homotopy excision it
is $(2n-c)$--cartesian for some (other) constant~$c$.  The second square
is cartesian since $\phi \simeq \check\phi$ is excisive.  By hypothesis
the maps in the upper left, upper right and lower left hand corners are
$(2n-c)$--connected, for some constant~$c$.  It follows that the map in the
lower right hand corner is also $(2n-c)$--connected, for some constant~$c$.
Hence $\theta(W_3)$ is a meta-stable equivalence.
\qed
\enddemo

By construction, the map $\theta(W)$ is the identity when $W = H_+$ and
when $W = *$.  Hence it is a strict equivalence whenever $W$ is weakly
contractible.  It follows that $\theta(\Sigma^k H_+)$ is a meta-stable
equivalence for all $k\ge0$, by induction on $k$, using lemma~10.3 applied
to the pushout square:
$$
\xymatrix{
\Sigma^k H_+ \ar[r] \ar[d] & C\Sigma^k H_+ \ar[d] \\
C\Sigma^k H_+ \ar[r] & \Sigma^{k+1} H_+ 
}
$$
Likewise it follows that $\theta(W)$ is a meta-stable equivalence for
all finite objects $W$ in $R(*, H)$, by induction on the number of free
$H$--cells in $W$, using lemma~10.3 applied to the pushout square associated
to the attachment of a free $H$--cell.  If $\phi$ satisfies the colimit
axiom, then this implies that $\theta(W)$ is a stable equivalence for
any based, free $H$--simplicial set $W$.  This completes the proof of~10.1.
\qed
\enddemo

The following form of the Brown--Whitehead representability theorem is
perhaps more familiar, although we shall not use it directly.

\proclaim{Proposition 10.4}
Let $\phi \: R(*, H) \to R(*, G)$ be a bounded below, linear functor
that satisfies the colimit axiom.  There is a natural chain of weak
equivalences
$$
\Omega^\infty(\partial\phi \wedge_H W) \simeq \phi(W)
$$
for all $W$ in $R(*, H)$.
\endproclaim

\demo{Proof}
The chain consists of the weak equivalence
$$
\Omega^\infty(\partial\phi \wedge_H W) @>\simeq>>
\Omega^\infty\{n \mapsto \check\phi(\Sigma^n W)\}
$$
of proposition~10.1, the weak equivalence
$$
\Omega^\infty\{n \mapsto \phi(\Sigma^n W)\} @>\simeq>>
\Omega^\infty\{n \mapsto \check\phi(\Sigma^n W)
$$
of~(9.2), and the weak equivalence
$$
\phi(W) @>\simeq>> \Omega^\infty\{n \mapsto \phi(\Sigma^n W)\}
$$
that follows since $\phi$ is linear.
\qed
\enddemo

\head The chain rule \endhead

Let $\Phi \: R(X) \to R(Y)$ and $\Psi \: R(Y) \to R(Z)$ be bounded below,
excisive functors, with $Y = \Phi(X)$ and $Z = \Psi(Y)$, such that $(X, x)$,
$(Y, y)$ and $(Z, z)$ are based connected simplicial sets.  We form $\phi
= j \circ \Phi \circ i$ and $\psi = j \circ \Psi \circ i$, as in~(9.1).
$$
\xymatrix{
R(X) \ar[r]^{\Phi} & R(Y) \ar[r]^{\Psi} \ar[d]<-0.5ex>_j & R(Z) \ar[d]^j \\
R(*, \Omega_x(X)) \ar[u]^i \ar[r]^{\phi} & R(*, \Omega_y(Y))
\ar[u]<-0.5ex>_i \ar[r]^{\psi} & R(*, \Omega_z(Z)) 
}
$$
We let $\Xi = \Psi \circ \Phi$ and $\xi = j \circ \Psi \circ
\Phi \circ i$, so that $\partial\xi \simeq \partial^x_z(\Xi)(X) =
\partial^x_z(\Psi\circ\Phi)(X)$ is the equivariant Goodwillie derivative
of the composite functor.

\proclaim{Proposition 11.1}
Suppose that $\Psi$ satisfies the colimit axiom.  Then there is a natural
chain of stable equivalences of free left $\Omega_z(Z)$--spectra with
right $\Omega_x(X)$--action
$$
\partial^x_z(\Psi \circ \Phi)(X)
\simeq \partial^y_z\Psi(Y) \wedge_{\Omega_y(Y)} \partial^x_y\Phi(X)
\,.
$$
\endproclaim

\demo{Proof}
By proposition~8.1 there is a natural weak equivalence $\xi = j \circ \Psi
\circ \Phi \circ i \simeq j \circ \Psi \circ i \circ j \circ \Phi \circ i
= \psi \circ \phi$.  Hence there is also a natural weak equivalence of
simplicial functors $\check\xi \simeq \check\psi \circ \check\phi$, and a
strict equivalence of equivariant Goodwillie derivatives $\partial\xi
\simeq \partial(\psi \circ \phi)$.

For brevity, let $H = \Omega_x(X)$ and $G = \Omega_y(Y)$.

The smash product of spectra is so constructed that in order to produce
a stable equivalence of spectra $\partial \psi \wedge_{G}
\partial \phi @>>> \partial(\psi \circ \phi)$ it suffices to define a
stable equivalence of bi-spectra
$$
\{m, n \mapsto \partial\psi_m \wedge_{G} \partial\phi_n \}
@>>> \{m, n \mapsto \partial(\psi \circ \phi)_{m, n} \}
$$
where
$$
\align
\partial\psi_m \wedge_{G} \partial\phi_n &=
\check\psi(\Sigma^m G_+) \wedge_{G} \check\phi(\Sigma^n H_+) \\
\partial(\psi \circ \phi)_{m, n} &=
(\check\psi \circ \check\phi)(\Sigma^m \Sigma^n H_+) \,.
\endalign
$$

Such a stable equivalence is provided by the following composite:
$$
\check\psi(\Sigma^m G_+) \wedge_{G} \check\phi(\Sigma^n H_+)
@>a>> \check\psi(\Sigma^m \check\phi(\Sigma^n H_+))
@>b>> (\check\psi \circ \check\phi)(\Sigma^m \Sigma^n H_+)
\,.
\tag 11.2
$$
Here the first map $a$ is a case of~(10.2), which induces a stable
equivalence (as $m \to \infty$) by proposition~10.1, in view of
lemma~8.4.  The second map $b$ is $\check\psi$ applied to the map
$$
\Sigma^m \check\phi(\Sigma^n H_+) @>>> \check\phi(\Sigma^m
\Sigma^n H_+) \,,
$$
which is~(10.2) applied to the case $W = \Sigma^m H_+$.  By proposition~10.1
again, this is a meta-stable equivalence (as $n \to \infty$).  Since
$\check\psi$ is bounded below, the second map $b$ is also a meta-stable
equivalence.
\qed
\enddemo

\proclaim{Theorem 11.3}
Let $\Phi \: \SS/X \to \SS$ and $\Psi \: \SS/Y \to \SS$ be bounded below,
stably excisive functors, with $Y = \Phi(X)$ and $Z = \Psi(Y)$, and suppose
that $\Psi$ satisfies the colimit axiom.  Write $Y = \coprod_{\beta \in B}
Y_\beta$ with each $Y_\beta$ connected, and choose base points $x \in X$,
$y_\beta \in Y_\beta$ and $z \in Z$.  Let $X_\alpha$ be the connected
component of $x$ in $X$, and let $Z_\gamma$ be the connected component of
$z$ in $Z$.  Then $\Psi \circ \Phi$ is also bounded below and
stably excisive, and there is a natural chain of stable equivalences
$$
\partial^x_z(\Psi \circ \Phi)(X)
\simeq \bigvee_{\beta \in B}
\partial^{y_\beta}_z\Psi(Y) \wedge_{\Omega_{y_\beta}(Y_\beta)}
\partial^x_{y_\beta}\Phi(X)
$$
of free left $\Omega_z(Z_\gamma)$--spectra with right
$\Omega_x(X_\alpha)$--action.
\endproclaim

\demo{Proof}
The composite $\Psi \circ \Phi$ is bounded below and stably excisive
by proposition~3.2.  By propositions~4.4 and~6.2 we can replace $\Phi$
and $\Psi$ by the bounded below, excisive functors $P'\Phi$ and $P'\Psi$,
respectively, without changing the derivatives of $\Phi$, $\Psi$ and
$\Psi \circ \Phi$ by more than a stable equivalence, and such that
$P'\Psi$ satisfies the colimit axiom.  Hence we can assume from the
beginning that $\Phi$ and $\Psi$ are excisive.  By lemma~7.1 and
proposition~7.3 there are stable equivalences
$$
\partial^x_z(\Psi \circ \Phi)(X)
\simeq \partial^x_z(\Psi_\gamma \circ \Phi^\alpha)(X_\alpha)
\simeq \bigvee_{\beta \in B} \partial^x_z(\Psi^\beta_\gamma
\circ \Phi^\alpha_\beta)(X_\alpha) \,.
$$
By proposition~11.1 the summand indexed by $\beta$ is stably equivalent to
$$
\partial^{y_\beta}_z (\Psi^\beta_\gamma)(Y_\beta)
\wedge_{\Omega_{y_\beta}(Y_\beta)} \partial^x_{y_\beta}
(\Phi^\alpha_\beta)(X_\alpha) \,.
$$
By lemma~7.1
each such term can be rewritten as $\partial^{y_\beta}_z \Psi(Y)
\wedge_{\Omega_{y_\beta}(Y_\beta)} \partial^x_{y_\beta}\Phi(X)$.
\qed
\enddemo

\proclaim{Theorem 11.4}
Let $E, F \: \SS \to \SS$ be bounded below, stably excisive functors,
with $Y = F(X)$ and $Z = E(Y)$, and suppose that $E$ satisfies
the colimit axiom.  Write $Y = \coprod_{\beta \in B}
Y_\beta$ with each $Y_\beta$ connected, and choose base points $x \in X$,
$y_\beta \in Y_\beta$ and $z \in Z$.  Let $X_\alpha$ be the connected
component of $x$ in $X$, and let $Z_\gamma$ be the connected component of
$z$ in $Z$.  Then the composite $E \circ F$ is bounded below and
stably excisive, and there is a natural chain of stable equivalences
$$
\partial^x_z(E \circ F)(X)
\simeq \bigvee_{\beta \in B}
\partial^{y_\beta}_z E(Y) \wedge_{\Omega_{y_\beta}(Y_\beta)}
\partial^x_{y_\beta} F(X)
$$
of free left $\Omega_z(Z_\gamma)$--spectra with right
$\Omega_x(X_\alpha)$--action.
\endproclaim

\demo{Proof}
This is the special case of theorem~11.3 when $\Phi = F \circ u$ and
$\Psi = E \circ u$.
\qed
\enddemo

\head Topological spaces \endhead

Let $\UU$ be the category of compactly generated topological spaces.
The geometric realization functor $|-| \: \SS \to \UU$ is left adjoint
to the total singular simplicial set functor $\Sing \: \UU \to \SS$.
There is a natural weak equivalence $|\Sing(X)| \to X$.

Given a weak homotopy functor $f \: \UU \to \UU$ we get a weak homotopy
functor $F \: \SS \to \SS$ by setting $F(X) = \Sing(f(|X|))$.  Let $Y =
f(X)$, $Y^1 = \Sing(f(X))$ and $Y^2 = F(\Sing(X))$.  There are natural
weak equivalences $|Y^2| \to |Y^1| \to Y$.  (The superscripts are simply
labels, and do not mean powers or skeleta.)

Choose base points $x \in X$ and $y^2 \in Y^2$.  Let $x \in \Sing(X)$,
$y^1 \in Y^1$ and $y \in Y$ denote the corresponding base points, via
the maps just mentioned.

The Goodwillie derivative of $f$ at $X$ with respect to the base points
$x \in X$ and $y \in Y$ is the spectrum
$$
\partial^x_y f(X) = \{
n \mapsto \hofib_y(f(X \vee_x S^n) \to Y) \} \,.
\tag 12.1
$$
It receives a natural strict equivalence from the spectrum
$$
|\partial^x_{y^2} F(\Sing(X))| = \{
n \mapsto |\hofib_{y^2}(F(\Sing X \vee_x S^n) \to Y^2| \} \,.
$$
The latter spectrum has a free left action by $|\Omega_{y^2}(Y^2_\beta)|$,
and a right action by $|\Omega_x(\Sing(X_\alpha))|$, where $X_\alpha$ is
the path component of $x$ in $X$ and $Y^2_\beta$ is the path component
of $y^2$ in $Y^2$.

It will be convenient to forgo the condition that the left action is
free.  We thus apply the forgetful functor $\SP(G) \to \SP^G$ from free
$G$--spectra to spectra with $G$--action, for $G = \Omega_{y^2}(Y^2_\beta)$.
Hence we will consider $\partial^x_y f(X)$ up to strict equivalence
as a spectrum with left $|\Omega_{y^2}(Y^2_\beta)|$--action and right
$|\Omega_x(\Sing(X_\alpha))|$--action.  We emphasize that our weak
equivalences of spectra with $G$--action are simply $G$--equivariant
maps that are stable equivalences, so that no fixed-point information
is retained.

We can always recover a strictly equivalent free $G$--spectrum by smashing
with $EG_+$, where $EG$ is a free, contractible $G$--space.  For example we
may take $EG = |\tilde\Gamma(Y^2_\beta)|$ as the geometric realization of
the principal bundle introduced in \S8.  Thus if $\bold L$ is a spectrum
with right $G$--action and $\bold M$ is a free left $G$--spectrum, thought
of as a spectrum with left $G$--action, the stable homotopy type of the
$G$--orbit spectrum $\bold L \wedge_G \bold M$ can be recovered as the
homotopy orbit spectrum
$$
\bold L \wedge_G (EG_+ \wedge \bold M) = \bold L \wedge_{hG} \bold M \,.
\tag 12.2
$$
We now switch to this notation.

Suppose that $e, f \: \UU \to \UU$ are bounded below, stably excisive
functors and that $e$ satisfies the colimit axiom.  Define $F(X) =
\Sing(f(|X|))$ and $E(Y) = \Sing(e(|Y|))$.  Let $Y = f(X)$ and $Z =
e(Y)$.  Let $Y = \coprod_{\beta \in B} Y_\beta$, $Y^1 = \coprod_{\beta
\in B} Y^1_\beta$ and $Y^2 = \coprod_{\beta \in B} Y^2_\beta$ be
the decompositions into path components.  (These simplicial sets are
weakly equivalent, so the indexing sets $B$ are equal.)  Choose base
points $x \in X$, $y^2_\beta \in Y^2_\beta$ and $z \in E(Y^2) = (E
\circ F)(\Sing(X))$.  Let $x \in \Sing(X)$, $y^1_\beta \in Y^1_\beta$,
$y_\beta \in Y_\beta$ and $z \in Z$ denote their images under the natural
weak equivalences.

The chain rule~11.4 for $E \circ F$ at $\Sing(X)$ then asserts that
there is a stable equivalence
$$
|\partial^x_z (E \circ F)(\Sing(X))| \simeq
\bigvee_{\beta \in B}
|\partial^{y^2_\beta}_z E(Y^2)|
\wedge_{h|\Omega_{y^2_\beta}(Y^2_\beta)|}
|\partial^x_{y^2_\beta} F(\Sing(X))| \,.
$$
(Geometric realization commutes with homotopy orbits since bisimplicial
sets can be realized in two stages, or at once.)  There are natural
weak equivalences
$$
\align
|\partial^x_{y^2_\beta} F(\Sing(X))|
&@>\simeq>> \partial^x_{y_\beta} f(X) \\
|\partial^{y^2_\beta}_z E(Y^2)|
&@>\simeq>> \partial^{y^2_\beta}_z e(Y^2)
@>\simeq>> \partial^{y_\beta}_z e(Y) \\
|\partial^x_z (E \circ F)(\Sing(X))|
&@>\simeq>> \partial^x_z (e \circ f)(X) \\
|\Omega_{y^2_\beta}(Y^2_\beta)|
&@>\simeq>> |\Omega_{y^1_\beta}(Y^1_\beta)| \,.
\endalign
$$
Hence we can summarize:

\proclaim{Theorem 12.3}
Let $e, f \: \UU \to \UU$ be bounded below, stably excisive functors,
with $Y = f(X)$ and $Z = e(Y)$, and suppose that $e$ satisfies the
colimit axiom.  Let $Y = \coprod_{\beta \in B} Y_\beta$ with each
$Y_\beta$ path connected, and choose base points $x \in X$, $y_\beta
\in Y_\beta$ and $z \in Z$.  Let $X_\alpha$ be the path component of
$x$ in $X$, and let $Z_\gamma$ be the path component of $z$ in $Z$.
Also set $\Omega_{y_\beta}(Y_\beta) = |\Omega_{y_\beta}(\Sing(Y_\beta))|$.
Then the composite $e \circ f$ is bounded below and stably excisive,
and there is a stable equivalence
$$
\partial^x_z (e \circ f)(X) \simeq
\bigvee_{\beta \in B} \partial^{y_\beta}_z e(Y) \wedge_{h\Omega_{y_\beta}(Y_\beta)}
\partial^x_{y_\beta} f(X)
$$
of spectra with left $\Omega_z(Z_\gamma)$--action and right
$\Omega_x(X_\alpha)$--action.
\endproclaim

\remark{Remark 12.4}
To be precise, this formula needs to be interpreted in line with the
weak equivalences above.  In particular, it is not
$\Omega_{y_\beta}(Y_\beta) = |\Omega_{y^1_\beta}(Y^1_\beta)|$
that really acts, but the naturally weakly equivalent topological
group $|\Omega_{y^2_\beta}(Y^2_\beta)|$.  And the action is
not really on the spectra $\partial^{y_\beta}_z e(Y)$ and
$\partial^x_{y_\beta} f(X)$, but on the naturally weakly equivalent
spectra $|\partial^{y^2_\beta}_z E(Y^2)|$ and
$|\partial^x_{y^2_\beta} F(\Sing(X))|$.  However, all the
constructions involved are weak homotopy invariant, so none of these
adjustments have any homotopy-theoretic significance.
\endremark

\eject
\head Examples \endhead

\example{Example 13.1}
Let $\id \: \UU \to \UU$ be the identity functor.  It is clearly
bounded below, and is stably excisive by homotopy excision.  For a
path connected space $X$, choose non-degenerate base points $x \in X$
and $y \in \id(X) = X$.  Let $P_y(X) = \{ \gamma \: I \to X \mid \gamma(0)
= y\}$ and $P^x_y(X) = \{\gamma \: I \to X \mid \gamma(0) = y, \gamma(1)
= x\}$.  There is a natural map
$$
\multline
\hofib_y(X \vee_x S^n \to X) = P_y(X) \cup_{P^x_y(X)} (P^x_y(X) \times
S^n) \\
@>>> P^x_y(X)_+ \wedge S^n = \Sigma^n P^x_y(X)_+
\endmultline
$$
which is a homotopy equivalence.  Hence there is a stable equivalence
$$
\partial^x_y(\id)(X) \simeq \Sigma^\infty P^x_y(X)_+ \,,
$$
with the natural left $\Omega_y(X)$--action and right $\Omega_x(X)$--action
given by composition of paths.
\endexample

\example{Example 13.2}
Let $K$ be a finite CW complex and consider the mapping space functor
$f(X) = X^K = \Map(K, X)$.  Then $f$ satisfies $E_1(d, \kappa)$ and
$E_2(2d+1, \kappa)$ (by homotopy excision) for all $\kappa$, where $d =
\dim(K)$.  Choose a non-degenerate base point $x \in X$, and base $Y =
X^K$ at the constant map $y$ taking $K$ to $\{x\}$.  There are homotopy
equivalences
$$
\hofib_y((X \vee_x S^n)^K \to X^K) \simeq \hofib_x(X \vee_x S^n \to X)^K
\simeq \Map(K, \Sigma^n \Omega_x(X)_+)
$$
and thus a stable equivalence
$$
\partial^x_y(X^K) \simeq \Map(K,
\Sigma^\infty \Omega_x(X)_+) \,.
$$
When $X^K$ is path connected, the group $\Omega_y(X^K) = \Map(K,
\Omega_x(X))$ acts from the left by pointwise multiplication, while
$\Omega_x(X)$ acts uniformly from the right.
\endexample

\example{Example 13.3}
Let $e(Y) = Q(Y_+) = \colim_n \Omega^n(\Sigma^n Y_+)$ be the (unreduced)
stable homotopy functor.  This functor is bounded below and excisive,
and satisfies the colimit axiom.  Choose a base point $y \in Y$, and
take any base point $z$ of $Z = Q(Y_+)$.  The pushout of $S^n \leftarrow
(Y \vee_y S^n)_+ \to Y_+$ is $*$, so there is a natural cartesian square
$$
\xymatrix{
Q((Y \vee_y S^n)_+) \ar[r] \ar[d] & Q(S^n) \ar[d] \\
Q(Y_+) \ar[r] & Q(*) 
}
$$
and a natural weak equivalence $\partial^y_z Q(Y_+)_n \to Q(S^n)$.
Thus there is a stable equivalence
$$
\partial^y_z Q(Y_+) \simeq \bold S \,,
$$
where $\bold S = \{n \mapsto Q(S^n)\}$ is the {\sl sphere spectrum}.
The left action of $\Omega_z(Q(Y_+))$ on $\partial^y_z Q(Y_+)_n$
pulls back from $*$ in the cartesian square above.  Likewise the right
action of $\Omega_y(Y)$ pulls back from the trivial action on $Q(S^n)$.
Hence both of these actions are trivial, up to homotopy.
\endexample

\example{Example 13.4}
Let $(e \circ f)(X) = Q(X^K_+)$ be the composite functor.
By the chain rule~12.3, its derivative at $X$ is
$$
\align
\partial^x_z Q(X^K_+) &\simeq \partial^y_z Q(Y_+) \wedge_{\Omega_y(Y)}
\partial^x_y(X^K) \\
&\simeq \bold S \wedge_{\Map(K, \Omega_x(X))} \Map(K, \Sigma^\infty
\Omega_x(X)_+) \\
&\simeq \Map(K, \Sigma^\infty \Omega_x(X)_+)_{h\Map(K, \Omega_x(X))} \,.
\endalign
$$
This assumes that $X^K$ is path connected.  The derivative of this
functor $X \mapsto Q(X^K_+)$ was first computed as the spectrum of
stable sections in a suitable Serre fibration in [\Cite{Go90}, \S2].
In the paper \cite{Kl95} the first author shows that the two descriptions
of this derivative are indeed equivalent.
\endexample

\enddocument